\documentclass[11pt]{article}
\usepackage{amsmath}
\usepackage{graphicx}
\usepackage{amsfonts}
\usepackage{amssymb}

\usepackage{amsfonts,latexsym,amstext}

\newtheorem{theorem}{Theorem}[section]
\newtheorem{lemma}[theorem]{Lemma}
\newtheorem{proposition}[theorem]{Proposition}
\newtheorem{definition}{Definition}[section]

\newtheorem{hypothesis}[theorem]{Hypothesis}

\newtheorem{remark}[theorem]{Remark}
\newtheorem{corollary}[theorem]{Corollary}

\makeatletter
\@addtoreset{equation}{section}

\makeatother

\def\qed{{\hfill\hbox{\enspace${ \square}$}} \smallskip}
\def\sqr#1#2{{\vcenter{\vbox{\hrule height .#2pt \hbox{\vrule
 width .#2pt height#1pt \kern#1pt \vrule
width .#2pt} \hrule height .#2pt}}}}
\def\square{\mathchoice\sqr54\sqr54\sqr{4.1}3\sqr{3.5}3}

\def\ds{\begin{displaystyle}}
\def\eds{\end{displaystyle}}
\def\dis{\displaystyle }
\def\<{\langle }
\def\>{\rangle }

\def\R{\mathbb R}

\def\E{\mathbb E}
\def\P{\mathbb P}

\def\calb{{\cal B}}

\def\cale{{\cal E}}
\def\calf{{\cal F}}
\def\calg{{\cal G}}

\def\calp{{\cal P}}
\def\calu{{\cal U}}

\def\call{{\cal L}}
\def\cals{{\cal S}}

\title{Reflected BSDEs and optimal control and stopping for infinite-dimensional systems }

\date{}
\author{Marco Fuhrman\\
Politecnico di Milano,
Dipartimento di Matematica\\
piazza Leonardo da Vinci 32, 20133 Milano, Italy\\
e-mail: marco.fuhrman@polimi.it\\
\\
Federica Masiero, Gianmario Tessitore\\
Dipartimento di Matematica e Applicazioni, Universit\`a di Milano Bicocca\\
via Cozzi 55, 20125 Milano, Italy\\
e-mail: federica.masiero@unimib.it,
gianmario.tessitore@unimib.it}

\begin{document}

\maketitle
\begin{abstract}
We introduce the notion of mild supersolution for an obstacle problem in an infinite dimensional Hilbert space. The minimal supersolution of this problem
is given in terms of a reflected BSDEs in an infinite dimensional Markovian framework. The results are applied to an optimal control and stopping problem.
\end{abstract}

\section{Introduction} \label{Introduction}

The connection between backward stochastic differential equations in $\mathbb{R}^n$ and semilinear parabolic PDEs is known since the seminal paper of Pardoux and Peng \cite{PaPe}. This result was extended to the case of reflected BSDEs and correspondingly of obstacle problem for PDEs in \cite{kakapapequ}.
Moreover it is also well known that the above equations are related to optimal stochastic control problems (in the first case)
and optimal stopping or optimal control/stopping problems in the second see \cite{YongZhou}. 
We notice that in the finite dimensional framework the above mentioned partial differential equations are intended either in classical sense (see \cite{PaPe}) or, more frequently,
in viscosity sense. 

On the other hand the relation between backward stochastic differential equations in infinite dimensional spaces,
optimal control of Hilbert valued stochastic 
evolution equations and parabolic equation on infinite dimensional spaces was investigated in \cite{FuTe} and in several successive papers.
In the above mentioned literature it appears that the concept of solution of the PDE 
has to be modified in the infinite dimensional case. Namely classical solutions require too much regularity while the theory
of viscosity solutions can be applied only in special cases with trace class noise and very regular value function
(see \cite{kesw}). The type of definition that was seen to fit the infinite dimensional framework 
and the BSDE approach is the classical notion of \textit{mild
solution}. 
Namely if we consider a semilinear parabolic PDE such as 
$$
\left\lbrace
\begin{array}{l}
\frac{\partial u}{\partial t}(t,x)=\call_{t}u\left(  t,x\right)
+\psi\left(  t,x,u\left(  t,x\right)  ,\nabla u\left(  t,x\right)\right)\\
\text{ \ \ \ \ }\qquad\qquad\qquad \qquad\qquad\qquad\qquad\qquad\qquad\qquad\qquad\qquad
t\in\left[  0,T\right]
,\text{ }x\in H\\
u(T,x)=\phi\left(  x\right)  ,
\end{array}
\right.  
$$
and $(P_{s,t})_{0\leq s \leq t\leq T}$ is the transition semigroup related to the second order differential operators
$(\call_{t})_{t\in [0,T]}$ then a function $u: [0,T]\times H \rightarrow \mathbb{R}$ is called a mild solution of the above PDE
whenever $u$ admits a gradient (in a suitable sense) and it holds:
$$
 u(s,x)= P_{s,t}[ u(t,\cdot)](x)
  +\int_s^t P_{s,\tau}\Big[
\psi (\tau, \cdot, u(\tau,\cdot),
\nabla u(\tau,\cdot))\Big](x) \; d\tau.
$$
 Large amount of literature has then extended the BSDE approach to control problems to several different situations both in the finite and in the infinite framework
but, at our best knowledge, the problem of relating reflected BSDEs in infinite dimensional spaces and obstacle problems for PDEs with infinitely many variables was never investigated.
The point is that it is not obvious how one should include the reflection term (which is not absolutely continuous with respect to Lebesgue measure on $[0,T]$)
into the definition of mild solution.

In this paper, inspired by the work of A. Bensoussan see \cite{Be}, to overcome such a difficulty, we propose the notion of \textit{mild-supersolution} (see Definition \ref{def-supersol}). To be more specific, our main result will be to prove  that if $(X^{s,x},Y^{s,x},Z^{s,x},K^{s,x})$
is the solution of the following forward backward system with reflected BSDE:
$$ \left\{\begin{array}{l}
dX^{s,x}_t  = AX^{s,x}_t +F(t,X^{s,x}_t)dt +G(t,X^{s,x}_t)dW_t
 \text{ \ \ \ }t\in\left[  s,T\right] \\
X^{s,x}_s  =x,\\
 - \dis dY^{s,x}_t=\psi(t,X^{s,x}_t,Y^{s,x}_t,Z^{s,x}_t)\;dt
+dK^{s,x}_t- Z^{s,x}_t \;dW_t,\qquad t\in [0,T],\\
  Y^{s,x}_T=\phi(X_T^{s,x}),\\
 Y^{s,x}_t\geq h(X_t^{s,x}),\\
 \int_0^T (Y^{s,x}_t-h(X_t^{s,x}))dK^{s,x}_t=0.
\end{array}\right.
$$
 setting  $u(t,x):= Y^{t,x}_t$ then
 $u$ is the minimal mild supersolution of the obstacle problem 
 \begin{equation}
\left\lbrace
\begin{array}{l}
\min\left(u(t,x)-h(x),
-\frac{\partial u}{\partial t}(t,x)-\call_{t}u\left(  t,x\right)
-\psi\left(  t,x,u\left(  t,x\right)  ,\nabla u\left(  t,x\right)G\left(  t,x\right)  \right) \right)\geq 0\\
\text{ \ \ \ \ }\qquad\qquad\qquad \qquad\qquad\qquad\qquad\qquad\qquad\qquad\qquad\qquad
t\in\left[  0,T\right]
,\text{ }x\in H\\
u(T,x)=\phi\left(  x\right)  ,
\end{array}
\right.  %
\end{equation}

Another issue that is considered in this paper is that we do not assume any nondegeneracy on the coefficient $G$ (and consequently 
any strong ellipticity on the second order differential operator in the PDE). Therefore we can not expect to have regular solutions of the obstacle problem. Thus we have to precise how the directional gradient $\nabla u G$ has to be intended. We choose here to employ the definition of generalized gradient (in probabilistic sense) introduced in \cite{FuTeGen}). It was proved in \cite{FuTeGen} that such generalized gradient  exists for all locally Lipschitz functions. In Theorem \ref{teo lip rifl} we prove that our candidate solution $u(t,x):= Y^{t,x}_t$ is indeed locally Lipschitz).
Moreover we notice that we work under general growth assumptions with respect to $x$ on the nonlinear term $\psi$ and on the final datum $\phi$.  This forces us to obtain $L^p$ estimates on the solution on the reflected BSDE that extend the ones proved in \cite{kakapapequ}.

The structure of the paper is the following. In section 2 we study reflected BSDEs obtaining the desired $L^p$ estimates and the local lipscitzianity with respect to the initial datum in the markovian framework. In section 3 we introduce the notion of minimal mild supersolution of the obstacle problem in the sense of the generalized gradient and we show how it is related to the reflected BSDEs. Finally in section 4 we apply the above results to an optimal control and stopping problem.
\section{Reflected BSDEs}

In a complete probability space $\left(  \Omega,\mathcal{F},\mathbb{P}\right)
$ we consider a cylindrical Wiener process $\left\{  W_{\tau},\tau
\geq0\right\}  $ in a Hilbert space $\Xi$ and $\left(  \mathcal{F}_{\tau}\right)  _{\tau\geq
0}$ is its natural filtration, augmented in the usual way.
We consider the following reflected backward stochastic
differential equation (RBSDE in the following):
\begin{equation}\label{RBSDEf}\left\{\begin{array}{l}
  \dis dY_t=-f(t,Y_t,Z_t)\;dt-dK_t+ Z_t \;dW_t,\qquad t\in [0,T],\\
  Y_T=\xi,\\
 Y_t\geq S_t,\\
 \int_0^T (Y_t-S_t)dK_t=0
\end{array}\right.
\end{equation}
for the unknown adapted processes $Y$, $Z$ and $K$. $Y$ and $K$ are real processes,
and $Z$ is a $\Xi^*$-valued process.
$Y$ and $Z$ are square integrable processes,
$Y$ admits a continuous modification and $K$ is a
continuous non-decreasing process with $K_0=0$.
The equation is understood in the usual integral way, namely:
\begin{equation}\label{RBSDEf-mild}
Y_t+\int_t^T Z_r \;dW_r=
\xi+\int_t^T f(r,Y_r,Z_r)\;dr
+K_T -K_t,\qquad t\in [0,T], \; \P-\hbox{a.s.}.
 \end{equation}
We also consider equation \ref{RBSDEf} with $f$ not depending on $(y,z)$:
\begin{equation}\label{RBSDEf-nodip}\left\{\begin{array}{l}
  \dis dY_t=-f(t)\;dt-dK_t+ Z_t \;dW_t,\qquad t\in [0,T],\\
  Y_T=\xi,\\
 Y_t\geq S_t,\\
 \int_0^T (Y_t-S_t)dK_t=0
\end{array}\right.
\end{equation}
In the following, if $E$ is a separable Hilbert space, $0<a<b$ and $p\geq 1$
 we denote by $L^p_\calp(\Omega\times[a,b], E )$ the space of
$E$-vauled
$\calf_t$-predictable processes $\ell$ s.t.
\[
\E \int_a^b \vert \ell(t)\vert^p dt <\infty.
\]
If $E=\R$ we write $L^p(\Omega\times[a,b] )$ instead of $L^p(\Omega\times[a,b], \R)$.

Moreover by $L^p_\calp(\Omega,C([a,b],E))$ we denote the subspace of $L^p_\calp(\Omega\times[a,b], E )$
given by processes admitting a continuous version and verifying
\[
\E \sup_{t\in [a,b]}\vert \ell(t)\vert^p <\infty.
\]
An analogous definition is given to $L^p_\calp(\Omega,C([a,b]))$

\noindent It proved in \cite[proposition 5.1]{kakapapequ}, that
if $f\in L^2_\calp(\Omega\times[0,T] )$, $\xi\in L^2(\Omega)$ and $\sup_{t\in[0,T]}S_t^+\in L^2(\Omega)$
then equation (\ref{RBSDEf-nodip}) admits a unique solution $(Y,Z,K)$
with $(Y,Z)\in  L^2_\calp(\Omega\times[0,T] )\times  L^2_\calp(\Omega\times[0,T],H )$, moreover
$Y$ admits a continuous version and
$\E\sup_{t\in [0,T]}\vert Y\vert^2<\infty$; finally $K_T\in L^2(\Omega)$.

In the following we need to prove regular dependence of the solution to the above equation with respect to parameters, namely the initial data of a related (forward) stochastic
differential equation. Due to the assumptions that we choose on the nonlinearity $\psi$ we will need $L^p$
estimates (both on the solution and on its approximations corresponding to suitable penalized approximating equations).

We make the following assumptions on the generator, on the final datum and on the obstable of the RBSDE (\ref{RBSDEf}):
\begin{hypothesis}\label{ip-f-YZ}
$f:(\Omega\times[0,T])\times \R\times\Xi\rightarrow \R$ is measurable with
respect to $\calp\times \calb\left(\R\times \Xi^*\right)$ ( where by $\calp$
we mean the predictable $\sigma$-algebra on $\Omega\times[0,T]$, and by $\calb(\Lambda)$
the Borel $\sigma$-algebra on any topological space $\Lambda$).

 Moreover $f$ is Lipschitz with respect to $y$ and $z$ uniformly in $t$ and $\omega$ and, for some $p\geq 2$
$$ \E\int_0^T\vert f(t,0,0)\vert^p < \infty
$$


The final data $\xi$ is $\calf_T$ measurable and $p$-integrable. 

Finally the obstacle $S$ is a continuous, $\calp$-meausurable, real valued process satisfying
 $$\E\sup_{t\in[0,T]}\vert S_t\vert^{2p-2}< \infty.$$
\end{hypothesis}

We notice that the integrability requests are not optimal (for instance we assume $p$-integrability jointly in $\Omega\times[0,T]$) for the generator $f$ and $2(p-1)$ integrability for the obstacle $S$). Nevertheless such assumptions are verified in the Markovian framework (see Section \ref{secMarkovianRBSDEs}) and will allow us to treat general obstacle problems under general assumptions (see Section \ref{sez-PDE}).

%
By a penalization procedure, we can prove the following:
\begin{theorem}\label{teo-stima-p-f}
 If hypothesis \ref{ip-f-YZ} holds true, equation (\ref{RBSDEf})
admits a unique adapted solution $(Y,Z,K)$ such that $Y$ admits a continuous version
and $K$ is non decreasing with $K_0=0$. Moreover $(Y,Z,K)$ satisfy
\begin{align}\label{stimaYZ_p-f}
\E&\sup_{t\in [0,T]}\vert Y_t\vert^p
+\E\left(\int_0^T\vert Z_t\vert^2 \,dt\right)^{p/2}+\E\vert K_T\vert^p\\
&\leq C\E\vert\xi\vert^p+
C\E\int_0^T\vert f(t,0,0)\vert^p\,dt+C\left(\E\sup_{t\in[0,T]}\vert S_t\vert^{2p-2}\right)^{p/(2p-2)}
. \nonumber
\end{align}
Where $C$ only depends on $T$ and on the Lipschitz constant of $f$.

\end{theorem}

We first need an analogous result on the corresponding penalized equation, that we now introduce.
Let us consider the following BSDE
\begin{equation}\label{penaliz-f}\left\{\begin{array}{l}
- \dis dY^{n}_t=f(t,Y^{n}_t,Z^{n}_t)\;dt
+n(Y^{n}_t-S_t))^-dt- Z^{n,s,x}_t \;dW_t,\qquad t\in [0,T],\\
Y^{n}_T=\xi.
\end{array}\right.
\end{equation}
It is shown in \cite{kakapapequ} that the penalized BSDE (\ref{penaliz-f}) admits a
unique solution $(Y^{n},Z^{n})$ in $L^2_\calp(\Omega,C([0,T]))\times L^2_\calp(\Omega\times[0,T],\Xi)$,
whose norm (in the above spaces) is uniformly bounded with respect to $n$. Moreover such a solution
$(Y^{n},Z^{n})$ converges in $L^2_\calp(\Omega,C([0,T]))\times L^2_\calp(\Omega\times[0,T],\Xi)$
to $Y^{s,x},Z^{s,x}$, solution of the RBSDE.
Next we want to prove an $L^p$-estimate, uniform with respect to $n$.

\begin{proposition}\label{prop-p-penalized-f}
 If hypothesis \ref{ip-f-YZ} holds true then equation (\ref{penaliz}) admits
admits a unique adapted solution $(Y^{n},Z^{n})$ such that $Y^{n}$ admits a continuous version.

\noindent Moreover $(Y,Z,K)$ satisfy
\begin{align}\label{stimaYZ_p-f-appr}
\E&\sup_{t\in [0,T]}\vert Y^n_t\vert^p
+\E\left(\int_0^T\vert Z^n_t\vert^2 \,dt\right)^{p/2}\\
& \leq C\E\vert\xi\vert^p+
C\E\int_0^T\vert f(t,0,0)\vert^p\,dt+C\left(\E\sup_{t\in[0,T]}\vert S_t\vert^{2p-2}\right)^{p/(2p-2)}
\nonumber .
\end{align}
Finally if $K^n_t=n \int_0^t (Y^n_s-S_s)^{-}ds$ then $K^n$ is an adapted, continuous, non-decreasing
proces satisfying
\begin{align}\label{stimaK_p-f-appr}
\E&\vert K^n_T\vert^p
\leq C\E\vert\xi\vert^p+
C\E\int_0^T\vert f(t,0,0)\vert^p\,dt+C\left(\E\sup_{t\in[0,T]}\vert S_t\vert^{2p-2}\right)^{p/(2p-2)}
\end{align}
where $C$ only depends on $p$, $T$ and on the Lipschitz constant of $f$.

\end{proposition}
\textbf{Proof.}  First of all we notice that we can always reduce ourselves to the case in which
\begin{equation}\label{monotonicita-f}
 \dfrac{y}{\vert y\vert }f(t,y,z)\leq \vert f(t,0,0)\vert +\mu \vert y\vert+\lambda \vert z\vert
\qquad \text{ with   }\mu+\lambda^2\leq 0.
\end{equation}
 Indeed, setting $\tilde Y^n_t=e^{at}Y^n_t,\,
\tilde Z^n_t=e^{at}Z^n_t$, we get that $\left(\tilde Y^n,\tilde Z^n\right)$
satisfies
\begin{equation*}\left\{\begin{array}{l}
- \dis d\tilde Y^{n}_t=e^{at} f(t,e^{-at}\tilde Y^{n}_t,e^{-at}\tilde Z^{n}_t)\;dt-a\tilde Y^{n}_t dt
+n(\tilde Y^{n}_t- \tilde S_t))^-dt\\
\qquad\qquad- \tilde Z^{n}_t \;dW_t,\qquad t\in [0,T],\\
Y^{n}_T=\xi.
\end{array}\right.
\end{equation*}
So the generator is given by
\[
 \tilde f(t,y,z):=e^{at} f(t,e^{-at}y,e^{-at}z)-ay,
\]
so by choosing $a$ sufficiently large (depending only on the Lipscitz  constant of $f$)
 we can assume $\mu+\lambda^2\leq -1$.
From now on we assume that (\ref{monotonicita-f}) holds true and for simplicity
we omit the superscript $\sim$ where necessary.

Moreover by $c$ we shall denote a constant that depends only on the Lipschitz constant of $f$, $T$ and $p$
and by $c(\delta)$  a constant that depends, beside the above parameters,  on an auxiliary  constant $\delta>0$. Their value can change from line to line.
\medskip
We apply It\^o formula to $\vert Y^n_t\vert^p$, $s\leq t\leq T$ and we get,
\begin{align*}
 -d\vert Y^n_t\vert^p= & p\vert Y^n_t\vert^{p-1}\hat Y^n_tf(t,Y^{n}_t,Z^{n}_t) dt +pn\vert Y^n_t\vert ^{p-1}\hat Y^n_t(Y^n_t-S_t)^-dt\\ \nonumber
&-p\vert Y^n_t\vert^{p-1}\hat Y^n_tZ^n_tdW_t-\dfrac{p(p-1)}{2}\vert Y^n_t\vert^{p-2}\vert Z^n_t\vert^2 dt.
\end{align*}
where $\hat Y^n_t:=\dfrac{Y^n_t}{\vert Y^n_t\vert}.$
Integrating between $s$ and $T$, $0\leq s\leq t\leq  T$, we get
\begin{align*}
&\vert Y_s^n\vert^p+ \dfrac{p(p-1)}{2}\int_s^T \vert Y^n_t\vert^{p-2} \vert Z^n_t\vert^2\, dt\\ \nonumber
 &   =\vert\xi\vert^p +p\int_s^T \vert Y^n_t\vert^{p-1}\hat Y^{n}_t
 f(t,Y^{n}_t,Z^{n}_t) dt  +np\int_s^T\vert Y^n_t\vert^{p-1}\hat Y^{n}_t (Y^n_t-S_t)^-\,dt\\ \nonumber
& -p\int_s^T \vert Y^n_r\vert^{p-1}\hat Y^{n}_tZ^n_tdW_t \\ \nonumber
 & \leq \vert\xi\vert^p +p\int_s^T \vert Y^n_t\vert^{p-1}
 \vert f(t,0,0)\vert dt +p\mu\int_s^T \vert Y^n_t\vert^{p}dt
+p\lambda\int_s^T \vert Y^n_t\vert^{p-1}\vert Z^n_t\vert dt \\ \nonumber
&+np\int_s^T\vert S_t\vert^{p-1} (Y^n_t-S_t)^-\,dt
 -p\int_s^T \vert Y^n_r\vert^{p-1}\hat Y^{n}_tZ^n_tdW_t \\ \nonumber
& \leq \vert\xi\vert^p +c\int_t^T \vert f(t,0,0)\vert^p dt
+p \int_s^T \vert Y^n_t\vert^{p}
 +p\mu\int_s^T \vert Y^n_t\vert^{p}dt
+\dfrac{p \lambda^2}{(p-1)}\int_s^T \vert Y^n_t\vert^{p} dt\\ \nonumber
&+\dfrac{p(p-1)}{4}\int_s^T \vert Y^n_t\vert^{p-2}\vert Z^n_t\vert^2 dt
+\sup_{t\in[s,T]} \vert S_t\vert^{p-1}n\int_s^T (Y^n_t-S_t)^-\,dt\\ \nonumber
& -p\int_s^T \vert Y^n_r\vert^{p-1}\hat Y^{n}_tZ^n_tdW_t, \nonumber
\end{align*}
where we have applied Young inequality.
So recalling that by (\ref{monotonicita-f}) $\mu+\lambda^2\leq 0$ and also since
 $p\geq 2$,  we get
\begin{align*}
&\vert Y_s^n\vert^p+ \dfrac{p(p-1)}{4}\int_s^T \vert Y^n_t\vert^{p-2} \vert Z^n_t\vert^2\, dt\\ \nonumber
 & \leq \vert\xi\vert^p +c\int_t^T \vert f(t,0,0)\vert^p dt
\\ \nonumber
 &+\sup_{t\in[s,T]} \vert S_t\vert^{p-1}n\int_s^T (Y^n_t-S_t)^-\,dt - c\int_s^T \vert Y^n_r\vert^{p-1}\hat Y^{n}_tZ^n_tdW_t \\ \nonumber
\end{align*}
By the penalized BSDE (\ref{penaliz-f}) in integral form we deduce that
\begin{equation}\label{def-di Kn}
   \int_s^Tn(Y^{n}_t-S_t)^-dt=-\xi+Y^{n}_t-\int_s^T f(t,Y^{n}_t,Z^{n}_t)\;dt
+ \int_s^TZ^{n}_t \;dW_t,
\end{equation}
and so
\begin{align}\label{stima-p-pen-0-f}
&\vert Y_s^n\vert^p+ \dfrac{p(p-1)}{4}\int_s^T \vert Y^n_t\vert^{p-2} \vert Z^n_t\vert^2\, dt\\ \nonumber
 & \leq \vert\xi\vert^p +c \int_s^T \vert f(t,0,0)\vert^p dt
  -p\int_s^T \vert Y^n_r\vert^{p-1}\hat Y^{n}_tZ^n_tdW_t \\ \nonumber
 &+  \left[\sup_{t\in[s,T]}\vert S_t\vert^{p-1}\left(-\xi+Y^{n}_s
-\int_s^Tf(t,Y^{n}_t,Z^{n}_t)\;dt + \int_s^TZ^{n}_t \;dW_t\right)\right]\\ \nonumber
  & \leq \vert\xi\vert^p +c\int_s^T \vert f(t,0,0)\vert^p dt
  -p\int_s^T \vert Y^n_r\vert^{p-1}\hat Y^{n}_tZ^n_tdW_t+c \sup_{t\in[s,T]} \vert S_t\vert^{p}+\dfrac{1}{2}\vert Y^{n}_s\vert^p  \\ \nonumber
 &+\sup_{t\in[s,T]} \vert S_t\vert^{p-1}\left[\vert \xi\vert
+\vert\int_s^Tf(t,Y^{n}_t,Z^{n}_t)\;dt\vert+
\vert\int_s^TZ^{n}_t \;dW_t\vert\right]
 \nonumber
\end{align}
Now we recall that,  by the $L^p$-estimates on BSDEs, see e.g. \cite{FuTe}, $(Y^n,Z^n)\in L^p_\calp(\Omega, C([0,T]) )\times L^p_\calp(\Omega, L^2([0,T],\Xi ))$,
and so the It\^o integral $\int_s^T \vert Y^n_r\vert^{p-1}\hat Y^{n}_tZ^n_tdW_t$ has null expectation.
Computing expectation  in the above inequality
\begin{align}\label{stima-p-pen-1-f}
&\dfrac{1}{2}\E\vert Y_s^n\vert^p+ \dfrac{p(p-1)}{4}\E\int_s^T \vert Y^n_t\vert^{p-2} \vert Z^n_t\vert^2\, dt
 \leq \E\vert\xi\vert^p + c \E\int_s^T \vert f(t,0,0)\vert^p dt\\ \nonumber
 &+\E\sup_{t\in[s,T]} \vert S_t\vert^{p-1}\left[\vert \xi\vert
+\vert\int_s^T f(t,Y^{n}_t,Z^{n}_t)\;dt\vert+
\vert\int_s^TZ^{n}_t \;dW_t\vert\right]\\ \nonumber
 &\leq C\E\vert\xi\vert^p +C\E\int_s^T \vert f(t,0,0)\vert^p dt +\left(\E\sup_{t\in[s,T]} \vert S_t\vert^{2(p-1)}\right)^{\frac{1}{2}}*\\ \nonumber
 &*\left(\E\left[ \vert\xi\vert+
 \int_s^T\left(|f(t,0,0)|+c\vert Y^n_t\vert +c\vert Z^n_t\vert\right)\;dt+\vert\int_s^TZ^{n}_t \;dW_t\vert\right]^2\right)^{\frac{1}{2}}
\end{align}
As already mentioned, it is well known that the penalized BSDE admits a unique solution whose norm is uniformly bounded in $L^2_\calp(\Omega,C([0,T]))\times L^2_\calp(\Omega\times[0,T],\Xi)$. Namely estimates in section 6 of \cite{kakapapequ} reed:
\[
 \E\sup_{s\in[0,T]}\vert Y_s^n\vert^2+\E\int_0^T \vert  Z^n_t\vert^2\, dt\leq c\E\vert\xi\vert^2 +c\E\int_0^T \vert f(t,0,0)\vert^2 dt,
\]
So plugging the above in
(\ref{stima-p-pen-1-f}) we get, also by the BDG inequality,
\begin{align}\label{stima-p-pen-2-f}
&\dfrac{1}{2}\E\vert Y_s^n\vert^p+ \dfrac{p(p-1)}{4}\E\int_s^T \vert Y^n_t\vert^{p-2} \vert Z^n_t\vert^2\, dt
 \leq \E\vert\xi\vert^p +c \E\int_0^T \vert f(t,0,0)\vert^p dt\\ \nonumber
 &+ c \left(\E\sup_{t\in[s,T]} \vert S_t\vert^{2(p-1)}\right)^{\frac{1}{2}}\E\left[\vert \xi\vert^2+\vert\int_0^T\vert f(t,0,0)\vert^2\;dt\vert+
  \int_s^T\vert Y^n_t\vert^2 +\int_s^T\vert Z^n_t\vert^2\;dt\right]^{\frac{1}{2}}\\ \nonumber
 &\leq c\E\vert\xi\vert^p +c\E\int_0^T \vert f(t,0,0)\vert^p dt
+c\left(\E\sup_{t\in[0,T]} \vert S_t\vert^{2(p-1)}\right)^{\frac{1}{2}}*\left(\E\vert\xi\vert^2 +C\E\int_0^T \vert f(t,0,0)\vert^2 dt \right)^{\frac{1}{2}}\nonumber
\end{align}
So we can deduce that
\begin{align}\label{stima-pen-intstoc-f}
 &\dfrac{p(p-1)}{4}\E\int_s^T \vert Y^n_t\vert^{p-2} \vert Z^n_t\vert^2\, dt\\
 &\leq  C\E\vert\xi\vert^p +C\E\int_0^T \vert f(t,0,0)\vert^p dt
+C\left(\E\sup_{t\in[0,T]} \vert S_t\vert^{2(p-1)}\right)^{\frac{1}{2}}*\left(\E\vert\xi\vert^2 +C\E\int_0^T \vert f(t,0,0)\vert^2 dt \right)^{\frac{1}{2}} \nonumber
\end{align}
By (\ref{stima-p-pen-0-f}), with $r$ in the place of $s$, such that $0\leq s\leq r\leq T$
we get
\begin{align*}
\vert Y_r^n\vert^p
  & \leq 2 \vert\xi\vert^p +c\int_0^T \vert f(t,0,0)\vert^p dt
  \\ \nonumber
  &  -2 p\int_r^T \vert Y^n_r\vert^{p-1}\hat Y^{n}_tZ^n_tdW_t+c \sup_{r\in[s,T]} \vert S_t\vert^{p}+  \\ \nonumber
 &+2 \sup_{r\in[s,T]} \vert S_r\vert^{p-1}\left[\vert \xi\vert+\vert\int_r^T f(t,Y^{n}_t,Z^{n}_t)\;dt\vert
+\vert\int_r^TZ^{n}_t \;dW_t\vert\right].
 \\ \nonumber
\end{align*}
By taking the supremum over the time $r$, by taking expectation
and with calculations in part similar to the ones we have
performed in (\ref{stima-p-pen-2-f})
we arrive at
\begin{align}\label{stima-p-pen-supY-f}
&\E\sup_{r\in [s,T]}\vert Y_r^n\vert^p\\ \nonumber
  & \leq 2\E\vert\xi\vert^p +c \E \int_s^T \vert f(t,0,0)\vert^p dt
  \\ \nonumber
  &  +c \E\sup_{r\in [s,T]}\vert\int_r^T \vert Y^n_t\vert^{p-1}\hat Y^{n}_tZ^n_tdW_t\vert+ c \E\sup_{t\in[s,T]} \vert S_t\vert^{p}\\ \nonumber
 &+2\left(\E\sup_{t\in[s,T]} \vert S_t\vert^{2(p-1)}\right)^{\frac{1}{2}}\left(\E\left[\vert \xi\vert^2+\int_s^T\vert f(t,Y^{n}_t,Z^{n}_t)\vert^2\;dt+
\sup_{t\in[s,T]}\vert\int_t^T Z^{n}_t \;dW_t\vert^2
\right]\right)^{\frac{1}{2}}
 \\ \nonumber
&\leq c\E\vert\xi\vert^p +c\E\int_s^T \vert f(t,0,0)\vert^p dt + c \E\sup_{r\in [s,T]}\int_r^T \vert Y^n_t\vert^{p-1}\hat Y^{n}_tZ^n_tdW_t\\ \nonumber
 \\ \nonumber
&+c \left(\E\sup_{t\in[s,T]} \vert S_t\vert^{2(p-1)}\right)^{\frac{1}{2}}\left(\E\left[\vert \xi\vert^2+\int_s^T\vert f(t,0,0)\vert^2\;dt\right]\right)^{\frac{1}{2}}
 \\
&\leq c\E\vert\xi\vert^p +c\int_s^T \vert f(t,0,0)\vert^p dt +c \E\left(\int_r^T \vert Y^n_t\vert^{2(p-1)}\vert Z^n_t\vert^2dt\right)^{\frac{1}{2}}  \nonumber
  \\ \nonumber
&+c \left(\E\sup_{t\in[s,T]} \vert S_t\vert^{2(p-1)}\right)^{\frac{1}{2}}\left(\E\left[\vert \xi\vert^2+\int_s^T\vert f(t,0,0)\vert^2\;dt\right]\right)^{\frac{1}{2}}   \nonumber
\\ \nonumber
&\leq c\E\vert\xi\vert^p +c\int_s^T \vert f(t,0,0)\vert^p dt +c \E\left(\sup_{t\in[s,T]}\vert Y^n_t\vert^{p}\int_s^T \vert Y^n_t\vert^{p-2}\vert Z^n_t\vert^2dt\right)^{\frac{1}{2}}
  \\ \nonumber
&+c \left(\E\sup_{t\in[s,T]} \vert S_t\vert^{2(p-1)}\right)^{\frac{1}{2}}\left(\E\left[\vert \xi\vert^2+\int_s^T\vert f(t,0,0)\vert^2\;dt\right]\right)^{\frac{1}{2}}
\\ \nonumber
&\leq c\E\vert\xi\vert^p + c \E \int_s^T \vert f(t,0,0)\vert^p dt +\dfrac{1}{2}\E \sup_{r\in [s,T]}\vert Y_r^n\vert^p + c\E \int_s^T \vert Y^n_t\vert^{p-2}\vert Z^n_t\vert^2dt
  \\ \nonumber
&+c \left(\E\sup_{t\in[s,T]} \vert S_t\vert^{2(p-1)}\right)^{\frac{1}{2}}\left(\E\left[\vert \xi\vert^2+\int_s^T\vert f(t,0,0)\vert^2\;dt\vert\right]\right)^{\frac{1}{2}}
  \nonumber
\end{align}
So we get, also by applying estimate (\ref{stima-pen-intstoc-f})
\begin{align}\label{stima-p-pen-supY-1-f}
\E\sup_{r\in [s,T]}\vert Y_r^n\vert^p
&\leq c\E\vert\xi\vert^p +c\E \int_s^T \vert f(t,0,0)\vert^p dt
  \\ \nonumber
 &+\left(\E\sup_{t\in[s,T]} \vert S_t\vert^{2(p-1)}\right)^{\frac{1}{2}}\left(\E\left[\vert \xi\vert^2+\int_s^T\vert f(t,0,0)\vert^2\;dt\right]\right)^{\frac{1}{2}}
\end{align}
Next  we estimate
$\E\left(\dis\int_s^T \vert Z^n_t\vert^2\, dt\right)^{\frac{p}{2}}$; we apply It\^o formula to $\vert Y^n_t\vert^2$, $s\leq t\leq T$ obtaining
\begin{equation*}
 d\vert Y^n_t\vert^2= -2Y^n_t f(t,Y^{n}_t,Z^{n}_t) dt -n Y^n_t(Y^n_t-S_t)^-dt
+2Y^n_tZ^n_tdW_t+\vert Z^n_t\vert^2 dt.
\end{equation*}
We integrate on $[s,T]$ and we raise to the power $\frac{p}{2}$:
\begin{align*}
&\vert Y_s^n\vert^p+ \left(\int_s^T \vert Z^n_t\vert^2\, dt\right)^{\frac{p}{2}}\\ \nonumber
 &\leq \vert\xi\vert^p +\left(2\int_s^T Y^n_t f(t,Y^{n}_t,Z^{n}_t) dt\right)^{\frac{p}{2}}  +\left(n\int_s^TY^n_t (Y^n_t-S_t)^-\,dt\right)^{\frac{p}{2}}\\ \nonumber
 &+\left(\int_s^T Y^n_rZ^n_tdW_t\right)^{\frac{p}{2}} \\ \nonumber
 &\leq \vert\xi\vert^p +\left(2\int_s^T Y^n_t f(t,Y^{n}_t,Z^{n}_t) dt\right)^{\frac{p}{2}}  +\left(\sup_{t\in[s,T]}\vert S_t\vert n\int_s^TY^n_t (Y^n_t-h(X_t))^-\,dt\right)^{\frac{p}{2}}\\ \nonumber
 &+\left(\int_s^T Y^n_rZ^n_tdW_t\right)^{\frac{p}{2}} \\ \nonumber
\end{align*}
Using the expression (\ref{def-di Kn}) for $n\dis\int_s^T (Y^n_t-S_t)^-\,dt$
that comes from the penalized BSDE (\ref{penaliz-f}), we get
\begin{align*}
& \left(\int_s^T \vert Z^n_t\vert^2\, dt\right)^{\frac{p}{2}}
  \\ \nonumber
  &\leq \vert\xi\vert^p +2\left|2\int_s^T Y^n_t f(t,Y^{n}_t,Z^{n}_t) dt\right|^{\frac{p}{2}}
 +\left|\int_s^T Y^n_tZ^n_tdW_t\right|^{\frac{p}{2}} \\ \nonumber &+\sup_{t\in[s,T]}\vert S_t\vert^{\frac{p}{2}} \left|-\xi+Y^{n}_t-\int_s^T f(t,Y^{n}_t,Z^{n}_t)\;dt
+ \int_s^TZ^{n}_t \;dW_t\right|^{\frac{p}{2}}\\ \nonumber
 &\leq \vert\xi\vert^p +\left(2\int_s^T \left(|Y^n_t| \,|f(t,0,0)|+\mu\vert Y^n_t\vert^2+
 \lambda\vert Y^n_t\vert \vert Z^n_t\vert\right) dt\right)^{\frac{p}{2}}
 +\left|\int_s^T Y^n_rZ^n_tdW_t\right|^{\frac{p}{2}} \\ \nonumber &
 +\sup_{t\in[s,T]}\vert S_t\vert^{\frac{p}{2}} \left(|\xi|+|Y^{n}_t|+\int_s^T\left(\vert f(t,0,0)\vert+
 \mu\vert Y^{n}_t\vert+\lambda\vert Z^{n}_t\vert\right)\;dt+ \left|\int_s^TZ^{n}_t \;dW_t\right|\right)^{\frac{p}{2}},\\ \nonumber
\end{align*}
Computing expectation, by BDG and Young inequalities, and by using
estimate (\ref{stima-p-pen-supY-1-f}), we get
\begin{align*}
& \E\left(\int_s^T \vert Z^n_t\vert^2\, dt\right)^{\frac{p}{2}}\\ \nonumber
&\leq \E\vert \xi\vert^p +c\E\sup_{r\in [s,T]}\vert Y^n_r\vert^p+
c\E\left(\int_s^T \vert f(t,0,0)\vert^2\,dt\right)^{\frac{p}{2}}
\\ \nonumber
 &+ \dfrac{1}{4} \left(\E\int_s^T\vert Z^n_t\vert^2 dt\right)^{\frac{p}{2}}
 +\E\left(\int_s^T \vert Y^n_tZ^n_t\vert^2dt\right)^{\frac{p}{4}}+ c \E\sup_{t\in[s,T]}\vert S_t)\vert^p
\\ \nonumber
 &\leq c \E\vert \xi\vert^p +c\E\int_0^T\vert f(t,0,0)\vert^p\,dt
+\dfrac{1}{2} \E\left(\int_s^T\vert Z^n_t\vert^2 dt\right)^{\frac{p}{2}}
 +c\E\sup_{r\in [s,T]}\vert Y^n_r\vert^p
 \end{align*}
 Concluding by estimate (\ref{stima-p-pen-supY-1-f}), we obtain:
\begin{align*}
\E&\left(\int_s^T \vert Z^n_t\vert^2\, dt\right)^{\frac{p}{2}}\leq c\E\vert \xi\vert^p +c\E\int_0^T\vert f(t,0,0)\vert^p\,dt \\ \nonumber
&+c\left(\E\sup_{t\in[s,T]} \vert S_t\vert^{2(p-1)}\right)^{\frac{1}{2}}\left(\E\left[\vert \xi\vert^2+\vert\int_s^T\vert f(t,0,0)\vert^2\;dt\vert\right]\right)^{\frac{1}{2}}
\end{align*}
and this concludes the estimate of $\E \left(\int_s^T \vert Z^n_t\vert^2\, dt\right)^{\frac{p}{2}}$.

\medskip
\noindent The estimate of $\E |K^n_T|^p $ is then easy consequence
 of the previous ones and of relation (\ref{def-di Kn}).
\qed

We are now ready to prove Theorem \ref{teo-stima-p-f}.

\medskip
\textbf{Proof of Theorem \ref{teo-stima-p-f}. } By \cite{kakapapequ}, section 6, we know that
$Y^n_t\uparrow Y_t$ and $$\E (\sup_{t\in [0,T]}(Y_t -Y^n_t))^2\rightarrow 0.$$ Thus choosing a suitable subsequence we can assume  the $\P$-a.s. convergence of $\sup_{t\in [0,T]}(Y_t -Y^n_t)$ towards 0. Consequently by Fatou Lemma and (\ref{stimaYZ_p-f-appr}) we get
$$
  \E\sup_{t\in [0,T]}\vert Y_t\vert^p \leq C\E\vert\xi\vert^p+
C\E\int_0^T\vert f(t,0,0)\vert^p\,dt+C\left(\E\sup_{t\in[0,T]}\vert S_t\vert^{2p-2}\right)^{p/(2p-2)}
 .$$
For what concerns the convergence of $Z^n$, again by
\cite{kakapapequ}, section 6, we already know that $Z^n\rightarrow Z$ in $L^2_\calp(\Omega\times [0,T])$,
and by proposition \ref{prop-p-penalized-f} we know that $Z^n$ is bounded in
$L^p_\calp(\Omega\times [0,T])$, so, extracting, if needed, a subsequence, we can assume that
such that $(Z^{n})$ converges weakly in $L^p_\calp(\Omega\times [0,T])$ and consequently also
weakly in $L^2_\calp(\Omega\times [0,T])$. Therefore the weak limit of $(Z^{n})$ in
$L^p_\calp(\Omega\times [0,T])$ must coincide with the strong limit $Z$ in $L^2_\calp(\Omega\times [0,T])$ topology. Consequently again by (\ref{stimaYZ_p-f-appr}) we have that $Z$ satisfies
$$
\E\left(\int_0^T\vert Z_t\vert^2 \,dt\right)^{p/2} \leq C\E\vert\xi\vert^p+
C\E\int_0^T\vert f(t,0,0)\vert^p\,dt+C\left(\E\sup_{t\in[0,T]}\vert S_t\vert^{2p-2}\right)^{p/(2p-2)}
 .$$
For what concerns $K$, by \cite{kakapapequ} we already know (see again \cite{kakapapequ}, section 6) that
$
 \E\vert K^{n}_T-K_T\vert^2\rightarrow 0,
$. The claim follows as before by Fatou lemma by extracting a subsequence that converges $\P$-a.s. and exploiting estimate (\ref{stimaYZ_p-f-appr}).

\qed

\subsection{Reflected BSDEs in a Markovian framework}
\label{secMarkovianRBSDEs}

Now we consider a RBSDE depending on a forward equation with values in
another real and separable Hilbert space $H$.
Namely, we consider the forward backward system
\begin{equation}\label{FRBSDE} \left\{\begin{array}{l}
dX^{s,x}_t  = AX^{s,x}_t +F(t,X^{s,x}_t)dt +G(t,X^{s,x}_t)dW_t
 \text{ \ \ \ }t\in\left[  s,T\right] \\
X^{s,x}_s  =x,\\
 - \dis dY^{s,x}_t=\psi(t,X^{s,x}_t,Y^{s,x}_t,Z^{s,x}_t)\;dt
+dK^{s,x}_t- Z^{s,x}_t \;dW_t,\qquad t\in [0,T],\\
  Y^{s,x}_T=\phi(X_T^{s,x}),\\
 Y^{s,x}_t\geq h(X_t^{s,x}),\\
 \int_0^T (Y^{s,x}_t-h(X_t^{s,x}))dK^{s,x}_t=0.
\end{array}\right.
\end{equation}
We denote the solution of the RBSDE in the above equation by $(Y^{s,x},Z^{s,x},K^{s,x})$, to stress
 the dependence on the initial conditions, or by $(Y,Z,K)$ if no confusion is
possible.

On the coefficients of the forward equation we make the following assumptions:
\begin{hypothesis}\label{ip-forward}
\begin{enumerate}
\item $A$ is the generator of a strongly continuous semigroup of linear operators $(e^{tA})_{t\geq 0}$;
\item  The mapping $F:\left[  0,T\right]  \times H\rightarrow H$ is measurable
and satisfies, for some constant $C>0$ and $0\leq\gamma<1$,
\begin{equation}%
\begin{array}
[c]{l}%
\left|  e^{sA}F\left(  \tau,x\right)  \right|  \leq Cs^{-\gamma}\left(
1+\left|  x\right|  \right)  ,\text{ \ \ \ }t\in\left[  0,T\right]  ,\\
\\
\left|  e^{sA}F\left(  \tau,x\right)  -e^{sA}F\left(  \tau,y\right)  \right|
\leq Cs^{-\gamma}\left|  x-y\right|  ,\text{ \ \ \ }s>0,\text{ }t\in\left[
0,T\right]  ,\text{ \ }x,y\in H.
\end{array}
\label{lipF}%
\end{equation}

\item $G$ is a mapping $\left[  0,T\right]  \times H\rightarrow L\left(
\Xi,H\right)  $ such that for every $v\in\Xi$, the map $Gv:\left[  0,T\right]
\times H\rightarrow H$ is measurable and for every $s>0,$ $\tau\in\left[
0,T\right]  $ and $x\in H$ we have $e^{sA}G\left(  \tau,x\right)  \in
L_{2}\left(  \Xi,H\right)  $. Moreover there exists $0<\theta<\frac{1}{2}$
such that
\begin{equation}%
\begin{array}
[c]{l}%
\left|  e^{sA}G\left(  \tau,x\right)  \right|  _{L_{2}\left(  \Xi,H\right)
}\leq Ls^{-\theta}\left(  1+\left|  x\right|  \right)  ,\\
\\
\left|  e^{sA}G\left(  \tau,x\right)  -e^{sA}G\left(  \tau,y\right)  \right|
_{L_{2}\left(  \Xi,H\right)  }\leq Ls^{-\theta}\left|  x-y\right|  ,\text{
\ \ \ }s>0,\text{ }\tau\in\left[  0,T\right]  ,\text{ \ }x,y\in H.
\end{array}
\label{condG}%
\end{equation}
\end{enumerate}
\end{hypothesis}
The next existence and uniqueness Proposition is proved in \cite{FuTe}.
\begin{proposition}
\label{esistenza sde}Under hypothesis \ref{ip-forward}, the forward equation in (\ref{FRBSDE})
admits a unique continuous mild solution. Moreover $\mathbb{E}\sup_{t
\in\left[  s,T\right]  }\left|  X_t^{s,x}  \right|  ^{p}%
<C_{p}\left(  1+\left|  x\right|  \right)  ^{p}$, for every $p\in\left(
0,\infty\right)  $, and some constant $C_{p}>0$.
\end{proposition}

We will work under the following assumptions on $\psi$:
\begin{hypothesis}\label{ip su psi}
The function $\psi:\left[  0,T\right]  \times H\times
\mathbb{R}\times\Xi\rightarrow\mathbb{R}$ is Borel measurable and
satisfies the following:

\begin{enumerate}
\item  there exists a constant $L>0$ such that
\[
\left|  \psi\left(  t,x,y_{1},z_{1}\right)  -\psi\left(  t,x,y_{2}%
,z_{2}\right)  \right|  \leq L\left(  \left|  y_{1}-y_{2}\right|  +\left|
z_{1}-z_{2}\right|  _{\Xi}\right)  ,
\]
for every $t\in\left[  0,T\right]  ,$ $x\in H,$ $y_{1},y_{2}\in\mathbb{R}$,
$z_{1},z_{2}\in\Xi$;

\item  for every $t\in\left[  0,T\right]  $, $\psi\left(  t,\cdot,\cdot
,\cdot\right)  $ is continuous $H\times\mathbb{R}\times\Xi^{\ast}%
\rightarrow\mathbb{R}$;

\item  there exists $L^{\prime}>0$ and $m\geq0$ such that
\[
\left|  \psi\left(  t,x_1,y,z\right) -\psi\left(  t,x_2,y,z\right) \right|  \leq L^{\prime}
\vert x_1-x_2\vert\left(1+\vert x_1\vert ^m +\vert x_2\vert^m+ |y|^m\right)\left(  1  +\left|  z\right|  _{\Xi}\right)  ,
\]
for every $t\in\left[  0,T\right]  $, $x_1,x_2\in H,$ $y\in\mathbb{R}$, $z\in
\Xi$.
\item as far as the final datum $\phi$ and the obstacle
$h$ are concerned there exists $L>0$ such that:
\begin{align*}
 &\vert \phi(x_1)-\phi(x_2)\vert \leq L
\vert x_1-x_2\vert(1+\vert x_1\vert ^m +\vert x_2\vert^m),\\
&\vert h(x_1)-h(x_2)\vert \leq L
\vert x_1-x_2\vert(1+\vert x_1\vert ^m +\vert x_2\vert^m),
\end{align*}
for all $x_1,x_2\in H$.
\end{enumerate}
\end{hypothesis}
We notice that hypothesis \ref{ip su psi} implies that, for all $p>0$
\begin{equation}
\label{ip su psi cons}
\left|  \psi\left(t,x,y,z\right) \right|  \leq L
\left(1+\vert x\vert ^{m+1} + \left|
y\right|  +\left|  z\right|  _{\Xi^{\ast}}\right) , \quad \vert \phi(x)\vert \leq L(1+\vert x\vert ^{m+1}),\quad
\vert h(x)\vert \leq L
(1+\vert x\vert ^{m+1}),
\end{equation}
for all $t\in[0,T]$, $x\in H$, $y\in\R$ $z\in\Xi$, and for all $p\geq 2$.
\begin{proposition}\label{prop-p-RBSDE}
Let hypotheses \ref{ip-forward} and \ref{ip su psi} hold true and fix $s\in [0,T]$, $x\in H$.
Then
%
%
 the RBSDE in (\ref{FRBSDE}) admits a unique adapted solution $(Y^{s,x},Z^{s,x},K^{s,x})$. Moreover
$Y^{s,x}$ admits a continuous version, $(K^{s,x})$ is continuous and  non-decreasing ($K^{s,x}_0=0$)
and, for all $p \geq 2$ there exists $C_p>0$ such that
\begin{equation}\label{stimaYZ_p}
\E\sup_{t\in [0,T]}\vert Y^{s,x}_t\vert^p
+\E\left(\int_0^T\vert Z^{s,x}_t\vert^2 \,dt\right)^{p/2}+\E|K^{s,x}_T|^p<C(1+\vert x\vert^{p(m+1)}).
\end{equation}
We consider also the penalized version of the RBSDE in (\ref{FRBSDE}):
\begin{equation}\label{penaliz}\left\{\begin{array}{l}
- \dis dY^{n,s,x}_t=\psi(t,X^{s,x}_t,Y^{n,s,x}_t,Z^{n,s,x}_t)\;dt
+n(Y^{n,s,x}_t-h(X^{s,x}_t))^-dt- Z^{n,s,x}_t \;dW_t,\qquad t\in [0,T],\\
Y^{n,s,x}_T=\phi(X_T^{s,x}).
\end{array}\right.
\end{equation}
The same holds for the penalized equation with constant $C$ independent on $n$.
\end{proposition}
\textbf{Proof.} It suffices to notice that by setting
\[
 f(t,y,z):=\psi (t,X^{s,x}_t, y,z),\; S_t:=h(X^{s,x}_t),\;\xi:=\phi(X^{s,x}_T)
\]
for all $t\in[0,T]$, $y\in\R$, $z\in\Xi$ and with $X^{s,x}$ solution to the forward equation
in the FBSDE (\ref{FRBSDE}),
by (\ref{ip su psi cons}) $f,\, h,\, S$ satisfy hypothesis \ref{ip-f-YZ}, and in particular:
\begin{align}\label{ip-f-YZ-cons}
 &\E\int_0^T \vert f(t,0,0)\vert^p dt =
\E\left|  \psi\left(  t,X_t,0,0\right) \right|^p  \leq c
\left(1+\vert x\vert ^{p(m+1)}\right)\\
&\E\sup_{t\in[0,T]}\vert S_t\vert^{2(p-1)}=\E\sup_{t\in[0,T]}\vert h(X_t)\vert^{2(p-1)}\leq c\left(1+\vert x\vert ^{2(p-1)(m+1)}\right)\\ \nonumber
&\E\vert \xi\vert^p=\E\vert\vert \phi(X_t)\vert^p \leq
(1+\vert x\vert ^{p(m+1)}).\\ \nonumber
\end{align}
So we can apply Proposition \ref{prop-p-penalized-f} and Theorem \ref{teo-stima-p-f} to obtain the claim \qed
\medskip
\begin{remark}\label{remark_makov:YZ} Notice that $(Y^{s,x}_t,Z^{s,x}_t)$ is independent on $\mathcal{F}_s$ so, fixed $0\leq \tau \leq s\leq T$ the compositions 
$$Y^{s,X^{\tau,x}_s}_t; \qquad Z^{s,X^{\tau,x}_s}_t,\qquad t\in [s,T] $$ 
are well defined. Moreover by uniqueness of the solution to the forward equation in (\ref{FRBSDE}) we have $X^{\tau,x}_t= X^{s,X^{\tau,x}_s}_t$ and consequently
$$Y^{s,X^{\tau,x}_s}_t=Y^{\tau,x}_t	\; \P-\hbox{a.s.},\; \forall t\in [s,T]$$
$$
Z^{s,X^{\tau,x}_s}_t=Z^{\tau,x}_t\; \P-\hbox{a.s. for a.e.}\; t\in [s,T]$$
\end{remark}
The next theorem is devoted to the local Lipschitz continuity of $Y^{s,x}$ with respect
to $x$.
\begin{theorem}\label{teo lip rifl}
 Let hypotheses \ref{ip-forward} and \ref{ip su psi} hold true and let $(Y^{s,x},Z^{s,x},K^{s,x})$
be the unique solution of the the RBSDE in (\ref{FRBSDE}).
Then there exists a constant
$L>0$ such that, $\forall x_1,x_2 \in H$,
\begin{equation}
\label{stima rbsde lip}
\vert Y^{s,x_1}_s-Y^{s,x_2}_s\vert\leq
L\left( 1+\vert x_1\vert^{m(m+1)}+\vert x_2\vert^{m(m+1)}\right)\vert x_1-x_2\vert.
 \end{equation}
\end{theorem}
\textbf{Proof.} We start by considering the generator
$\psi$ differentiable, namely for every $t\in\left[  0,T\right]$
we assume that $\psi\left(  t,\cdot,\cdot
,\cdot\right)  \in\calg(H\times\mathbb{R}\times\Xi^{\ast},\mathbb{R})$. The idea is to prove that, in the case
of smooth (differentiable) coefficients, the solution of
the penalized equation (\ref{penaliz}) is differentiable
with respect to $x$, and the derivative is bounded
uniformly with respect to $n$ so that in particular
we get local lipschitz continuity of $Y^{n,s,x}_s$ with respect to $x$,
 that
 is preserved
as $n\rightarrow\infty$.

\noindent In order to work in a ``smooth'' framework,
in the penalized BSDE (\ref{penaliz}) instead of considering
the penalizing term $n(y-h)^-$, we have to consider a smooth penalizing term,
namely we consider a function $\gamma:\R\rightarrow\R$, such that $\gamma\in C^{\infty}_b(\R)$
\begin{align*}
 &\gamma(y)=0 \text{  for }y\geq0,\qquad \gamma(y)>0 \text{  for }y<0\\
&\gamma(y)=-y \text{  for }y\leq -1,
\qquad \dot{\gamma}(y)<0 \text{  for }y<0.
\end{align*}
Notice that to construct $\gamma$ it is enough to set $\gamma(y)=\int_0^{-y} \ell(r) dr$ with
\begin{align*}
 &\ell(r)=0 \text{  for }r\leq 0,\qquad \ell(r)>0 \text{  for }r>0,\qquad \ell(r)=1 \text{  for }r\geq 1
 ,\qquad \int_0^{1} \ell(r) dr=1.
\end{align*}

So we consider the following ``smooth'' penalized BSDE
\begin{equation}\label{penaliz-smooth}\left\{\begin{array}{l}
- \dis dY^{n,s,x}_t=\psi(t,X^{s,x}_t,Y^{n,s,x}_t,Z^{n,s,x}_t)\;dt
+n\gamma(Y^{n,s,x}_t-h(X^{s,x}_t))dt- Z^{n,s,x}_t \;dW_t,\qquad t\in [0,T],\\
Y^{n,s,x}_T=\phi(X_T^{s,x}),
\end{array}\right.
\end{equation}
and we notice that estimates obtained in proposition \ref{prop-p-RBSDE}
are still true for the pair of processes $(Y^{n,s,x},Z^{n,s,x})$
solution of equation (\ref{prop-p-RBSDE}).

\noindent Notice that it is still true that 
$|y|^{p-1} \hat {y}\gamma(y-s)\leq |s|^{p-1}\gamma(y-s) $ for all $y,s,\in \R$.

\noindent By \cite{FuTe} we know that we can differentiate
 $(Y^{n,s,x},Z^{n,s,x})$  with respect to $x$, and that $(\nabla_x Y^{n,s,x},\nabla_x Z^{n,s,x})$
 is the solution of the BSDE (to be intended in mild form):
\begin{equation*}\left\{\begin{array}{l}
- \dis d\nabla_xdY^{n,s,x}_t=\nabla_x\psi(t,X^{s,x}_t,Y^{n,s,x}_t,Z^{n,s,x}_t)\nabla_x X_t^{s,x}\;dt
+\nabla_y\psi(t,X^{s,x}_t,Y^{n,s,x}_t,Z^{n,s,x}_t)\nabla_x Y_t^{n,s,x}\;dt\\
\qquad\qquad\qquad+n\dot\gamma(Y^{n,s,x}_t-h(X^{s,x}_t))(\nabla_x Y_t^{n,s,x}-\nabla h(X^{s,x}_t)\nabla_xX^{s,x}_t)dt\\
\qquad\qquad\qquad+\nabla_z\psi(t,X^{s,x}_t,Y^{n,s,x}_t,Z^{n,s,x}_t)\nabla_x Z_t^{n,s,x}\;dt
-\nabla_x Z^{n,s,x}_t \;dW_t,\qquad t\in [s,T],\\
\nabla_x Y^{n,s,x}_T=\nabla \phi(X_T^{s,x})\nabla_x X_T^{s,x}.
\end{array}\right.
\end{equation*}
where (see again \cite{FuTe}) $\nabla_x X^{s,x}$ is the mild solution to the following
forward equation
\begin{equation*}\left\{\begin{array}{l}
 \dis d \nabla_xX^{s,x}_t=A\nabla_x X_t^{s,x}\;dt
+\nabla_x F(t,X^{s,x}_t)\nabla_x X_t^{s,x}\;dt
+\nabla_x G(t,X^{s,x}_t)\nabla_x X_t^{s,x} \;dW_t,\qquad t\in [s,T],\\
\nabla_x X^{s,x}_s=I,
\end{array}\right.
\end{equation*}
 $I:H\rightarrow H$ being the identity operator in $H$.

\noindent 
We set $\tilde{\P}:= \mathcal{E}_T \P$, with
\begin{equation}\label{Girsanov density}
\mathcal{E}_T=\exp\left(\int_s^T \nabla_z\psi(t,X^{s,x}_t,Y^{n,s,x}_t,Z^{n,s,x}_t)\,dW_t-\frac{1}{2}
 \int_s^T \vert\nabla_z\psi(t,X^{s,x}_t,Y^{n,s,x}_t,Z^{n,s,x}_t)\vert^2 dt  \right).
\end{equation} 
By the Girsanov theorem $\tilde{\P}$ is a probability measure
equivalent to the original one $\P$ (recall that  by hypothesis \ref{ip su psi}, $ \nabla_z$ is bounded)
and 
\[
\tilde{W}_\tau=-\int_s^\tau \nabla_z\psi(t,X^{s,x}_t,Y^{n,s,x}_t,Z^{n,s,x}_t) dt
+W_\tau,\quad  s\leq\tau\leq T
\]
is a  $\tilde{\P}$-cylindrical Wiener process.

In $(\Omega,\calf,\tilde\P)$ the pair $(\nabla_xY^{n,s,x},\nabla_xZ^{n,s,x})$ solve
the following BSDE for $t\in [s,T]$:
\begin{equation}\label{penalizdiffle-smooth}\left\{\begin{array}{l}
- \dis d\nabla_xdY^{n,s,x}_t=\nabla_x\psi(t,X^{s,x}_t,Y^{n,s,x}_t,Z^{n,s,x}_t)\nabla_x X_t^{s,x}\;dt\\
\qquad\qquad\qquad
+\nabla_y\psi(t,X^{s,x}_t,Y^{n,s,x}_t,Z^{n,s,x}_t)\nabla_x Y_t^{n,s,x}\;dt\\
\qquad\qquad\qquad
+n\dot\gamma(Y^{n,s,x}_t-h(X^{s,x}_t))(\nabla_x Y_t^{n,s,x}-\nabla h(X^{s,x}_t)\nabla_xX^{s,x}_t)dt
-\nabla_x Z^{n,s,x}_t \;d\tilde{W}_t, \\
\nabla_x Y^{n,s,x}_T=\nabla \phi(X_T^{s,x})\nabla_x X_T^{s,x},
\end{array}\right.
\end{equation}
Multiplying $\nabla_x Y^{n,s,x}_t$ by $\displaystyle \exp\left\lbrace\int_s^t
 (\nabla_y\psi(t,X^{s,x}_\sigma,Y^{n,s,x}_\sigma,Z^{n,s,x}_\sigma)+
n\dot\gamma(Y^{n,s,x}_\sigma-h(X^{s,x}_\sigma)))\;d\sigma\right\rbrace$
and writing the obtained equation in $t=s$ we get:
\begin{align}\label{nablaYn-explicit}
& \nabla_xY^{n,s,x}_s \nonumber \\
&=\E\bigg[\mathcal{E}_T\int_s^T \exp\left\lbrace\int_s^\tau
 \nabla_y\psi(t,X^{s,x}_\sigma,Y^{n,s,x}_\sigma,Z^{n,s,x}_\sigma)+
n\dot\gamma(Y^{n,s,x}_\sigma-h(X^{s,x}_\sigma))\;d\sigma\right\rbrace\\ \nonumber
&\quad\Big(\nabla_x\psi(\tau,X^{s,x}_\tau,Y^{n,s,x}_\tau,Z^{n,s,x}_\tau)\nabla_x X_t^{s,x}
-n\dot\gamma(Y^{n,s,x}_\tau-h(X^{s,x}_\tau))\nabla_xX^{s,x}_\tau\Big)\;d\tau\bigg]\\ \nonumber
&+\E\Big[\mathcal{E}_T \exp\left\lbrace\int_s^T
 \nabla_y\psi(t,X^{s,x}_\sigma,Y^{n,s,x}_\sigma,Z^{n,s,x}_\sigma)+
n\dot\gamma(Y^{n,s,x}_\sigma-h(X^{s,x}_\sigma)) \;d\sigma\right\rbrace
\nabla\phi(X_T^{s,x})\nabla_xX_T^{s,x}\Big],\\ \nonumber
\end{align}
so that,  since $\dot\gamma\leq0$ and
$\nabla_y\psi$ is bounded by hypothesis
\ref{ip su psi}, point 1,
\begin{align*}
&\vert\nabla_xY^{n,s,x}_s\vert\\ \nonumber
&\leq c \E\Big[\cale_T \vert\nabla\phi(X_T^{s,x})\nabla_xX_T^{s,x}+
\int_s^T\nabla_x\psi(\tau,X^{s,x}_\tau,Y^{n,s,x}_\tau,Z^{n,s,x}_\tau)\nabla_x X_\tau^{s,x}\,d\tau
\vert\Big]\\ \nonumber
&+c \E\Big[\cale_T \vert\int_s^T \exp\left\lbrace\int_s^\tau
 n\dot\gamma(Y^{n,s,x}_\sigma-h(X^{s,x}_\sigma))\;d\sigma\right\rbrace
\left(-n\dot\gamma(Y^{n,s,x}_\tau-h(X^{s,x}_\tau))\right)\nabla h(X^{s,x}_\tau)
\nabla_xX^{s,x}_\tau\;d\tau\vert\Big]\\ \nonumber
&=I+II,
\end{align*}
We start by estimating I.  Here and in the following we again denote by $c$ a constant whose value can vary from line to line and that may depend on $T$, on the coefficients
$A,\,F,\,G,\,\psi,\,,h,\,\phi$, on $p$ but not on $n$ and $x$. 
$$
I\leq c\E\Big[\cale_T\vert\nabla\phi(X_T^{s,x})\nabla_xX_T^{s,x}
\vert\Big]+c\E\Big[\cale_T\int_s^T\vert\nabla_x\psi(\tau,X^{s,x}_\tau,Y^{n,s,x}_\tau,Z^{n,s,x}_\tau)
\nabla_x X_\tau^{s,x} \vert\,d\tau\Big].
$$
Taking into account that $\E \mathcal{E}_T^p\leq c$,
by Holder inequality, with $p,\,q,\,r$ conjugate exponents  $p>1,\,1<q<2,\,qm>2$,
(where $m$ is the same as in hypothesis \ref{ip su psi}) we get:
$$ \E\Big[\cale_T\vert\nabla\phi(X_T^{s,x})\nabla_xX_T^{s,x}
\vert\Big]\\ \nonumber
 \leq c\Big(\E\Big[\vert\nabla\phi(X_T^{s,x})\vert^q\Big]\Big)^{1/q}
 \Big(\E\Big[\vert\nabla_xX_T^{s,x}
 \vert^r\Big]\Big)^{1/r} \leq c (1+\vert x\vert^m),
$$
where we  have used  the estimate on $\nabla_xX_T^{s,x}$
stated in \cite{FuTe}, proposition 3.3.

\noindent 
In a similar way we can estimate (for $q>2$)
\begin{align*}
  &\E\Big[\cale_T\int_s^T\vert\nabla_x\psi(\tau,X^{s,x}_\tau,Y^{n,s,x}_\tau,Z^{n,s,x}_\tau)
 \nabla_x X_\tau^{s,x} \vert\,d\tau\Big]\\ \nonumber
 &\leq c\E\Big[\cale_T\int_s^T\left(1+\vert X^{s,x}_\tau\vert^{m}+\vert Y^{n,s,x}_\tau\vert^{m}\right)
  \vert Z^{n,s,x}_\tau\vert
    \left(1+\vert\nabla_x X_\tau^{s,x}\vert\right)\,d\tau\Big]\\ \nonumber
  &\leq c\left(\E\Big[\left(1+\sup_{\tau\in[s,T]}\vert X^{s,x}_\tau\vert^{mq}+\sup_{\tau\in[s,T]}\vert Y^{n,s,x}_\tau\vert^{mq}\right)   \left(1+\sup_{\tau\in[s,T]}\vert\nabla_x X_\tau^{s,x}\vert^{q}\right)
 \left(\int_s^T\vert Z^{n,s,x}_\tau\vert
  \,d\tau\right)^{q}\Big]\right)^{1/q} 
  \\ \nonumber
 &\leq c \left(\E\Big[1+\sup_{\tau\in[s,T]}\vert X^{s,x}_\tau\vert^{2mq}
 + \sup_{\tau\in[s,T]}\vert Y^{n,s,x}_\tau\vert^{2mq}\Big]\right)^{1/2q}
 \left(\E\Big[\int_s^T\vert Z^{n,s,x}_\tau\vert
 \,d\tau\Big]^{2q}\right)^{1/2q} 
 \\ \nonumber
&\leq c\left(1+\vert x\vert^{m(m+1)}\right).
\end{align*}
where we have used estimates \ref{stimaYZ_p} and  Proposition \ref{esistenza sde}.

\noindent For what concerns $II$, let $p,q$ and $\bar p,\,\bar q$ be two pairs of conjugate exponents, and let $$l(\tau):= -n\dot\gamma(Y^{n,s,x}_\tau-h(X^{s,x}_\tau))\geq 0, \;\tau \in [s,T]$$
Then
\begin{align*}
&\E\Big[\cale_T\Big\vert\int_s^T \exp\Big(-\int_s^{\tau }l_\sigma\;d\sigma\Big )
\nabla h(X^{s,x}_\tau)\nabla_xX^{s,x}_\tau\;d\tau
\Big \vert\Big]\\ \nonumber
&\leq c 
\left(\E\Big[\Big(1+\sup_{\tau \in[s,T]}\vert X^{s,x}_\tau\vert^m\Big)
\sup_{\tau \in[s,T]}\vert\nabla_xX^{s,x}_\tau\vert\int_s^T\exp\Big(-\int_s^{\tau }l_\sigma\;d\sigma\Big ) l_{\tau} d\tau
\Big]^q\right)^{1/q}
\\ \nonumber
& \leq c 
\left(\E\Big[\left(1+\sup_{\tau \in[s,T]}\vert X^{s,x}_\tau\vert^{m\bar p q}\right)
\sup_{\tau \in[s,T]}\vert\nabla_xX^{s,x}_\tau\vert^{\bar p q}\Big]\right)^{1/(\bar p q)} \left(\E\Big[\int_s^T\exp\Big(-\int_s^{\tau }l_\sigma\;d\sigma\Big ) l_{\tau}
\;d\tau
\Big]^{q\bar q}\right)^{1/(q \bar q)}
\\ \nonumber
& \leq c 
\left(1+\vert x\vert^{m}\right)
\left(\E\Big[ 1- \exp\Big(-\int_s^{\tau }l_\sigma\;d\sigma\Big )
\Big]^{q\bar q}\right)^{1/(q\bar q)}  \leq c 
\left(1+\vert x\vert^{m}\right)
\end{align*}
where in the last passage we have used that
\[
 \int_s^T \exp\left\lbrace-\int_s^\tau l(\sigma)\;d\sigma\right\rbrace l(\tau)\,d\tau
=1-\exp\left\lbrace-\int_s^T l(\sigma)\;d\sigma\right\rbrace
\]
So 
\begin{equation} \label{stimanabla-n}
\vert\nabla_xY^{n,s,x}_s\vert\leq c\left(1+\vert x\vert^{m(m+1)}\right),
\end{equation}
where $c$ may depend on $T$, on the coefficients
$A,\,F\,G,\,\psi,\, h,\,\phi$, but not on $n$.
By (\ref{stimanabla-n}) we get that $\forall\,x,y\in H$
\begin{equation*}
\vert Y^{n,s,x}_s- Y^{n,s,y}_s\vert\leq c \vert x-y\vert(1+\vert x\vert^{m(m+1)}+\vert y\vert^{m(m+1)}).
\end{equation*}
By letting $n\rightarrow\infty$, arguing as in section 6 in \cite{kakapapequ}
finally get the desired Lipschitz continuity of $Y^{s,x}_s$:
\begin{equation}
\vert Y^{s,x}_s- Y^{s,y}_s\vert\leq c \vert x-y\vert(1+\vert x\vert^{m(m+1)}+\vert y\vert^{m(m+1)}),\qquad\forall\,x,\,y\,\in H.
\end{equation}
Finally we have to remove the assumption of differentiability
on the coefficient $\psi,\,h,\,\phi$ in the reflected BSDE.
Since $\psi $ is Lipschitz continuous with respect to $y$ and $z$,
and $\forall\,t\,,y\,,z\,\in[0,T]\times\R\times\Xi$, $\psi(t,\cdot,y,z)\,,h\,,\phi$
are locally Lipschitz continuous with respect to $x$, then by taking
their inf-sup convolution $(\psi_k,\,\phi_k,\,h_k)_{k\geq 1}$ we obtain differentiable functions where the derivative
is bounded by the Lipschitz constant in the Lipschitz case, and the derivative
has the polynomial growth imposed by the locally Lipschitz growth, see e.g. \cite{DP3}
for the notion of inf-sup convolution, and
\cite{Mas} and \cite{Mas1} for the use of inf-sup convolutions in the Lipschitz and locally Lipschitz case.
So in particular the growth of the derivatives the inf-sup convolutions is uniform
with respect to $k$: it follows that Lipschitz and locally Lipschitz constants are
uniform with respect to $k$, and this allows to pass to the limit as $k\rightarrow\infty$
and to preserve Lipschitz and locally Lipschitz properties.
Coming into more details, we denote by $(Y^{n,k,s,x},Z^{n,k,s,x},K^{n,k,s,x})$
the solution of the penalized RBSDEs with regularized coefficients:
\begin{equation}\label{penaliz-smooth-infsup}\left\{\begin{array}{l}
- \dis dY^{n,k,s,x}_t=\psi_k(t,X^{s,x}_t,Y^{n,k,s,x}_t,Z^{n,k,s,x}_t)\;dt
+n\gamma(Y^{n,k,s,x}_t-h_k(X^{s,x}_t))dt- Z^{n,k,s,x}_t \;dW_t,\\

\qquad\qquad\qquad\qquad\qquad\qquad\qquad\qquad\qquad\qquad\qquad\qquad\qquad\qquad\qquad\qquad\qquad t\in [0,T],\\
Y^{n,k,s,x}_T=\phi_k(X_T^{s,x}),
\end{array}\right.
\end{equation}
By the previous calculations we get that $\forall\,x,y\in H$
\begin{equation*}
\vert Y^{n,k,s,x}_s- Y^{n,k,s,y}_s\vert\leq c \vert x-y\vert(1+\vert x\vert^{m(m+1)}
+\vert y\vert^{m(m+1)}),
\end{equation*}
where $c$ does not depend on $n$ nor on $k$.
By standard results on BSDEs (see \cite{FuTe_Bismut}) we know that
\[
( Y^{n,k,s,x},  Z^{n,k,s,x})\rightarrow (Y^{n,s,x}, Z^{n,s,x})\;
\text{in}\; L^p_\calp(\Omega, C([0,T]))\times L^p_\calp(\Omega\times[0,T],\Xi),
\]
where $(Y^{n,s,x}, Z^{n,s,x})$
is solution to the smooth penalized BSDE (\ref{penaliz-smooth}). 
In particular $Y^{n,k,s,x}_s \rightarrow Y^{n,s,x}_s$. Finally proceeding as in \cite{kakapapequ}
if we let $n\rightarrow \infty$, we have already recalled
that  $ Y^{n,s,x}\uparrow Y^{s,x}$, since for this monototne convergence what matters is the monotonicity of the penalization term corrisponding to $K$, and
we get the desired Lipschitz continuity of $Y^{s,x}_s$:
\begin{equation}
\vert Y^{s,x}_s- Y^{s,y}_s\vert\leq c \vert x-y\vert(1+\vert x\vert^{m(m+1)}+\vert y\vert^{m(m+1)}),\qquad\forall\,x,\,y\,\in H.
\end{equation}
\qed.

\begin{remark}\label{remark-lim-RBSDE}
Notice that if $h$ and $\phi$ are bounded and lipschitz continuous functions,
if for every $s\in[0,T]$, $\sup_{x\in H}\vert\psi(s,x,0,0)\vert<\infty$ and as a function of $x$,
$\psi$ is lipschitz continuous uniformly with respect to the other variables, that is hypothesis
\ref{ip su psi}, point 3 holds true with $m=0$, then by repeating the same argument in proposition
\ref{prop-p-RBSDE}, we can prove that the processes $Y^{s,x},\,Z^{s,x}$
are bounded processes with respect to $x$, that is namely
\begin{equation}\label{stima-bounded-YZ}
 \E\sup_{t\in [0,T]}\vert Y^{s,x}_t\vert^p
+\E\left(\int_0^T\vert Z^{s,x}_t\vert^2 \,dt\right)^{p/2}<C.
\end{equation}
\end{remark}

\section{Obstacle problem for a semilinear parabolic PDE: solution via RBSDEs}
\label{sez-PDE}

In this section we consider an obstacle problem for a semilinear
PDE in an infinite dimensional Hilbert space $H$ and we
solve it in a suitable sense by means of reflected BSDEs.
An informal description is as follows: we study an obstacle problem
of the following form%
\begin{equation}
\left\lbrace
\begin{array}{l}
\min\left(u(t,x)-h(x),
-\frac{\partial u}{\partial t}(t,x)-\call_{t}u\left(  t,x\right)
-\psi\left(  t,x,u\left(  t,x\right)  ,\nabla u\left(  t,x\right)G\left(  t,x\right)  \right) \right)=0\\
\text{ \ \ \ \ }\qquad\qquad\qquad \qquad\qquad\qquad\qquad\qquad\qquad\qquad\qquad\qquad
t\in\left[  0,T\right]
,\text{ }x\in H\\
u(T,x)=\phi\left(  x\right)  ,
\end{array}
\right.  \label{obstacle}%
\end{equation}
where $G:\left[  0,T\right]  \times H\rightarrow L\left(  \Xi,H\right)  $,
and $\nabla u\left(  t,x\right)G(t,x ) $ is the directional generalized
gradient of $u$ with respect to $x$,
see \cite{FuTeGen}, section 3, and the following
for the definition of generalized gradient.
For a function $f:H\rightarrow \R$, the operator $\call_{t}$ is
formally defined by%
\[
\call_{t}f\left(  x\right)  =\frac{1}{2}Trace\left(  G\left(
t,x\right)  G^{\ast}\left(  t,x\right)  \nabla^{2}f\left(  x\right)  \right)
+\langle Ax,\nabla f\left(  x\right)  \rangle_{H}+\langle F\left(  t,x\right)
,\nabla f\left(  x\right)  \rangle_{H},
\]
and it arises as the generator of an appropriate Markov process $X$ in $H$.

\noindent More precisely if $X$ is the mild solution to the
stochastic differential equation in $H$%
\begin{equation}
\left\{
\begin{array}
[c]{ll}%
dX_t^{s,x}  =\left[  AX_t^{s,x}   +F\left( t,X_t^{s,x}   \right)  \right]  dt+G\left(t,X_t^{s,x} \right)
dW_t  , & t\in\left[ s,T\right] \\
X_s^{s,x}   =x, & x\in H,
\end{array}
\right.  \label{sdes}%
\end{equation}
where $T>0$ is fixed.
For $t\in[s,T]$ we denote by $P_{s,t}$ the transition semigroup 
\[
 P_{s,t}[\phi](x)=\E\phi(X_t^{s,x}).
\]
where  $\phi:H\rightarrow\R$ is  bounded and measurable.

\noindent Note that $\call_t$ is formally the generator of the transition semigroup
$(P_{s,t})_{t\in[s,T]}$.
This leads us to consider solutions of the obstacle problem (\ref{obstacle})
in mild sense, as we are going to state.

\subsection{The generalized directional gradient}
\label{subsez-gradgen}
We observe that, under our assumptions, it is reasonable to expect that function 
$u$ is locally  Lipschitz but not that it is differentiable.

\noindent To this aim, we briefly show an example where the value function of a deterministic optimal stopping problem is not differentiable. Let us consider, as state equation without control,
\begin{equation*}
\left\{
\begin{array}
[c]{l}%
dX_{t}^{ s,x}=0
\\
X_{s}^{ s,x}=x\in \R
\end{array}
\right.
\end{equation*}
We consider the following cost functional:
\begin{equation*}
J\left(  s,x,\tau, \right)  =\phi\left(  X_{T}^{u,s,x} \right)\chi_{\left\lbrace \tau =T\right\rbrace} +
 h\left(\tau , X_{\tau}^{u,s,x}\right) \chi_{\left\lbrace\tau <T\right\rbrace} ,
\end{equation*}
So the value function is given by
\begin{equation*}
u\left(  x \right)  =\sup_\tau\left(\phi\left(  X_{T}^{u,s,x} \right)
\chi_{\left\lbrace \tau =T\right\rbrace}  +
 h\left(\tau , X_{\tau}^{u,s,x}\right) \chi_{\left\lbrace\tau <T\right\rbrace} \right)
=\sup_{x\in\R}(\phi(x),h(x)),
\end{equation*}
and it is evident that, even if the data are differentiable,
the value function may fail to be differentiable.

 Notice that in the above example and statement we take into account that we allow degeneracy of the noise. The issue of differentiability of $u$ when noise is non degenerate is very interesting but falls out of the scope of the present work. 
 
To take into account the lack of regularity of $u$ the derivative $\nabla u$ must not appear in the precise formulation of the problem.
Indeed it will be substituted by the notion of generalized gradient, whose definition is given in the next subsection.

We start by giving the definition of generalized gradient

%

%
%

\begin{theorem}\label{ggrad}
  Assume that Hypothesis \ref{ip-forward} holds
and that $u:[0,T]\times H\to \R$ is a Borel
measurable function satisfying, for some $r>0$
 \begin{equation}\label{iposuu}
  |u(t,x)-u(t,x')|\leq c (1 + |x|+|x'|)^r |x-x'|.
  \end{equation}
    Then
   there exists a Borel measurable function $\zeta:
        [0,T]\times H \rightarrow \Xi^*$ with the following
properties.
        \begin{enumerate}
        \item[(i)]
For every
  $s\in [0, T]$, $x\in H$ and $p\in [2,\infty)$,
\begin{equation}\label{sommabilitadizeta}
  \E\int_s^T|\zeta(\tau,X_\tau^{t,x})|^p\;d\tau
        <+\infty.
\end{equation}

\item[(ii)] For $\xi\in\Xi$, $x\in H$   and
  $0\le s\le T'<T$ the processes
  $\{u(t,X_t^{s,x}),t\in[s,T]\}$
  and $W^\xi$ admit a joint quadratic
  variation on the interval $[s,T']$ and
        $$
 \langle u(\cdot,X_\cdot^{s,x}),W^\xi
\rangle_{[s,T']}=\int_s^{T'}
\zeta(t,X_t^{s,x})\xi\; dt,\qquad \P-a.s.
$$
\item[(iii)] Moreover there exists
a Borel measurable function
$\rho: [0,T]\times H \rightarrow H^*$
such that for all $t\in [s,T]$ and all $x\in H$
$$\zeta(t,X_t^{s,x})=\rho(t,X_t^{s,x})G(t,X_t^{s,x})
\quad \hbox{$\P$-a.s.
for a.a. $t\in [s,T]$}
$$
\end{enumerate}
\end{theorem}
\textbf{Proof.}  
The proof is given in \cite{FuTeGen}, section 4. In that paper
it is also noticed, see remark 3.1, that uniqueness can be stated in the following sense:
if $\hat{\zeta}$ is another function
with the stated properties then for $0\leq s\leq t \leq T$ and $x\in H$ we
have
$
\zeta(t,X_t^{s,x})= \hat{\zeta}(t,X_t^{s,x})$, 
$\P-{\rm a.s. \;for \; a.a.\;} t\in [s,T].
$
\qed

\begin{definition}\label{defggrad}
Let $u:[0,T]\times H\to \R$ be a
  Borel measurable function satisfying
(\ref{iposuu}). The family of all measurable functions $\zeta:
[0,T]\times H \rightarrow \Xi^*$ satisfying properties (i) and
(ii)  in Theorem \ref{ggrad} will be called the {\em generalized
directional gradient of $u$} and denoted by $\widetilde{\nabla}^G
u$.
\end{definition}

\subsection{Mild solutions of the obstacle problem in the sense of the generalized
directional gradient}
\label{subsez-mildsol}

Having defined the generalized directional gradient, we are in the position to
give the precise definition of supersolution for the problem (\ref{obstacle}).

\begin{definition}\label{def-supersol} We say that a Borel measurable function $\bar u:[0,T]\times H\to\R$ is
a mild supersolution of the obstacle problem (\ref{obstacle}) in the sense of
the generalized directional  gradient  if the following holds:\begin{enumerate}
    \item for some $C>0, r\geq 0$ and
   for every $s\in [0,T]$, $x,y\in H$
    $$
   |\bar u(s,x)-\bar u(s,y)|\leq C|x-y|(1+|x|+|y|)^{r},
    \qquad |u(s,0)|\leq C;
$$
\item for every $s\in [0,T]$, $x\in H$,
$$
\bar u(s,x)\geq h(x);
$$
    \item
    for all
     $0\leq s\leq t\leq T$ and $x\in H$
 \begin{equation}\label{defdisolvarq}
\bar u(s,x)\geq P_{s,t}[ u(t,\cdot)](x)
  +\int_s^t P_{s,\tau}\Big[
\psi (\tau, \cdot,\bar u(\tau,\cdot),
\zeta(\tau,\cdot))\Big](x) \; d\tau,
\end{equation}
where $\zeta$ is an arbitrary element of
the generalized gradient $\widetilde{\nabla}^G\bar u$;
\item
$
 \bar u(T,\cdot)=\phi. $

\end{enumerate}
\end{definition}

We are now ready to state the main result of this paper.

\begin{theorem}\label{teo-solmild} Assume that hypotheses \ref{ip-forward} and \ref{ip su psi} hold true.
Let us define
\begin{equation}\label{def-di-u}
 u(s,x)=Y_s^{s,x},
\end{equation}
where $(Y^{s,x},Z^{s,x})$ is solution to the reflected BSDE in (\ref{FRBSDE}). Then $u$ is a mild
supersolution in the sense of the generalized directional gradient
for the obstacle problem (\ref{obstacle}).

 Moreover $u$ is minimal in the sense that given any $\bar u$, supersolution of (\ref{obstacle}) in the sense
of definition \ref{def-supersol}, it holds
$u(s,x)\leq\bar u(s,x),$  and $ s\in [0,T], \;x\in H$

Finally, if in addition $\sup_{s\in[0,T],x\in H}\vert\psi(s,x,0,0)\vert<\infty$ and $\phi$ and $h$ are bounded
then $u$ is also bounded.
 \end{theorem}
\textbf{Proof.}  By theorem \ref{teo lip rifl} by defining $u(s,x):=Y_s^{s,x}$,
$u$ has the regularity required in definition
\ref{def-supersol}, point 1, and moreover points 2 and 4 immediately follow
since $Y$ is solution to the RBSDE in (\ref{FRBSDE}).

 For what concerns point 3 of definition
\ref{def-supersol}, since $Y$ is solution to the reflected BSDE, we get
\begin{equation}\label{u-mild-RBSDE}
 u(s,x)=Y_t^{s,x}
  +\int_s^t
\psi (\tau,X_\tau^{s,x},Y_\tau^{s,x},
Z_\tau^{s,x}) \; d\tau +K_t^{s,x}-K_s^{s,x}-\int_s^t Z_\tau^{s,x}\,dW_\tau,
\end{equation}
 Fixed $\xi\in\Xi$, let us consider the
joint quadratic variation of both sides of (\ref{u-mild-RBSDE})
with $W^\xi$. Proposition 2.1 in \cite{FuTeGen} and Theorem \ref{ggrad} yield that
$\widetilde{\nabla}^G u$  exists and letting $\zeta\in\widetilde{\nabla}^G u$, we have
\[
 \<u(\cdot, X_\cdot^{s,x}), W^\xi_\cdot\>_{[s,t]}=\int_s^t
\zeta(\sigma,X_\sigma^{s,x})\xi\; d\sigma,
\]
where
\[
 W^\xi_t:=\int_s^t
\<\xi,\; dW_\sigma\>, \qquad 0\leq s\leq t\leq T.
\]
On the other hand by the Markov property stated in Remark \ref{remark_makov:YZ}
\[
 u(t, X_t^{s,x})= Y^{t, X_t^{s,x}}_t=Y_t^{s,x}
\]
and since $Y$ is solution to the RBSDE in (\ref{FRBSDE}) we deduce:
\[
 \<Y^{s,x}_\cdot, W_\cdot\>_{[s,t]}=\int_s^t
Z_\sigma^{s,x}\xi\; d\sigma.
\]
So, by these two expression of the joint quadratic variation of $u(\cdot, X_\cdot^{s,x})$
and $W^\xi$ we get
\begin{equation}\label{identifquasi}
\int_s^t
\zeta(\sigma,X_\sigma^{s,x})\xi\; d\sigma=
\int_s^t Z_\sigma^{t,x}\xi\; d\sigma,
\end{equation}
$\P$-a.s. Since both sides of (\ref{identifquasi}) are
continuous with respect to $t$, it follows that, $\P$-a.s.,
they coincide for all $t\in [s,T]$. This
implies that
$ \zeta(t,X_t^{s,x})=Z_t^{s,x}$, $\P$-a.s. for a.a.
$t\in [s, T]$. 
Therefore equation (\ref{u-mild-RBSDE}) can be rewritten as 
\begin{equation}\label{u-mild-RBSDE-zeta}
 u(s,x)=Y_t^{s,x}
  +\int_s^t
\psi (\tau,X_\tau^{s,x},u(\tau,X_\tau^{s,x}),
\zeta(\tau,X_{\tau}^{s,x})) \; d\tau +K_t^{s,x}-K_s^{s,x}-\int_s^t Z_\tau^{s,x}\,dW_\tau,
\end{equation}
By taking the conditional expectation $\E^{\calf_s}$ and since $K$ is
a nondecreasing process, we get
\begin{equation}\label{u-mild-RBSDE-exp}
  u(s,x)\geq Y_t^{s,x}
  +\int_s^t P_{s,\tau}\Big[
\psi (\tau, \cdot, u(\tau,\cdot),
 \zeta(\tau,\cdot))\Big](x) \; d\tau,
\end{equation}
 and we have proved that $u$ is a mild supersolution along the Definition  \ref{def-supersol}

\noindent We have to prove that $u$ is the minimal supersolution.
Let $\bar u$ be any supersolution and let us define $\bar Y_t^{s,x}=\bar u(t,X_t^{s,x})$.
Then for every $\sigma\in[s,t]$, with $0\leq s\leq t$, by point 3 of definition \ref{def-supersol},
having replaced $x$ with $X_\sigma^{s,x}$ which is $\calf_\sigma$-measurable,
\begin{equation}\label{u-supersol-mild}
\bar u(\sigma,X_\sigma^{s,x})
\geq\E^{\calf_\sigma}\bar u(t,X_t^{\sigma,X_\sigma^{s,x}})
  +\E^{\calf_\sigma}\int_\sigma^t
\psi (\tau,X_\tau^{\sigma,X_\sigma^{s,x}},\bar Y_\tau^{\sigma,X_\sigma^{s,x}},
\bar\zeta(\tau,X_\tau^{\sigma,X_\sigma^{s,x}}) \; d\tau.
\end{equation}
So it turns out that
$$
\left(L_\sigma^{s,x}\right)_{\sigma\in[s,T]}:=\left(-\bar u(\sigma,X_\sigma^{s,x})-\int_s^\sigma \psi \left(\tau,X_\tau^{s,x},\bar Y_\tau^{s,x},
\bar\zeta(\tau,X_\tau^{s,x})\right) \; d\tau\right)_{\sigma\in[s,T]}.
$$
is a submartingale. By  hypothesis
\ref{ip su psi} on $\psi$,
by the growth property of
$u$ as required in definition \ref{def-supersol}, point 1, by relation \ref{sommabilitadizeta} and finally by Proposition \ref{esistenza sde}
we get that $L^{s,x}$ is a uniformly integrable martingale, so it is of class
(D) and the Doob-Meyer decomposition can be applied,
see e.g. Definition 4.8 and Theorem 4.10 in Chapter 1 of
\cite{KaSh}.
So $L^{s,x}$ can be decomposed into:
$$
L_\sigma^{s,x}=\bar M_\sigma^{s,x}+\bar K_\sigma^{s,x},
$$
where $\bar K^{s,x}$ is an integrable nondecreasing process such that $\bar K_s^{s,x}=0$, and
$\bar M^{s,x}$ is a uniformly integrable martingale.
Moreover, see \cite{DelMey}, Chapter VII relation (15.1), since
\[
 \E\sup_{\sigma\in[s,T]}\vert L_\sigma^{s,x}\vert^2<\infty
\]
we have
\[
 \bar K^{s,x}_T\in L^2(\Omega).
\]
Notice that we are working in a
complete probability space filtered with the filtration generated by the Wiener process, so
by the martingale representation theorem, see again \cite{KaSh} and \cite{DP1} for
its infinite dimensional version, there exists a process
$\bar Z\in L^2_\calp(\Omega\times[s,T];L_2(\Xi,\R))$
such that
\[
 \bar M_\sigma^{s,x}=-\Big[u(s,x)+\int_s^\sigma \bar Z_\tau^{s,x}\,dW_\tau\Big].
\]
We finally get $\forall \sigma\in[s,T]$
\begin{equation}\label{quasi-RBSDE}
 u(s,x)=\bar u(\sigma,X_\sigma^{s,x})+\int_s^\sigma
\psi (\tau,X_\tau^{s,x},\bar Y_\tau^{s,x},
\bar\zeta(\tau,X_\tau^{s,x}) \; d\tau+\bar K_\sigma^{s,x}-\bar K_s^{s,x}
-\int_s^t\bar Z_\tau^{s,x}\,dW_\tau,
\end{equation}
that is, $\forall\, 0\leq s\leq t\leq T$
\begin{equation}\label{quasi-RBSDE-1}
\bar{Y}_t^{s,x}=\bar{Y}_T^{s,x}
  +\int_t^T
\psi (\tau,X_\tau^{s,x},\bar Y_\tau^{s,x},
\bar\zeta(\tau,X_\tau^{s,x})) \; d\tau +\bar K_t^{s,x}-\bar K_s^{s,x}
-\int_s^t\bar Z_\tau^{s,x}\,dW_\tau.
\end{equation}
Finally we have to identify
$ \tau \zeta(\tau,X_\tau^{s,x})$ with $\bar Z_\tau^{s,x}$, $\P$-a.s. for a.a.
$\tau\in [s, T]$. To this aim, for $\xi\in\Xi$, let us consider the
joint quadratic variation of both sides of (\ref{quasi-RBSDE})
with $W^\xi$. Notice that the finite variation term $K$ does not
give any contribution to the joint quadratic variation with $W^\xi$;
so Proposition 2.1 in \cite{FuTeGen} and Theorem \ref{ggrad} yield, for $s\leq
\sigma <T$ and $\zeta\in\widetilde{\nabla}^G u$,
\begin{equation}\label{identifquasi-bis}
\int_s^\sigma\zeta(\tau,X_\tau^{s,x})\xi\; d\tau=
\int_s^\sigma \bar Z_\tau^{s,x}\xi\; d\tau,
\end{equation}
$\P$-a.s.. Since both sides of (\ref{identifquasi-bis}) are
continuous with respect to $\sigma$, it follows that, $\P$-a.s.,
they coincide for all $\sigma\in [s,T]$. This
implies that
$ \zeta(\sigma,X_\sigma^{s,x})=\bar Z_\sigma^{s,x}$, $\P$-a.s. for a.a.
$\sigma\in [s, T]$.
So we get that, by defining, $\bar Y_s^{s,x}:=\bar u(s,x)$ and
$\bar Y^{s,x}:=\bar u(\cdot, X^{s,x})$, the couple of processes $(\bar Y^{s,x},\bar Z^{s,x})$
solves the following problem
\begin{equation}\label{quasi-RBSDE-bis}\left\{\begin{array}{l}
 - \dis dY_t^{s,x}=\psi(t,X_t^{s,x},\bar Y_t^{s,x},\bar Z_t^{s,x})\;dt
+d\bar K_t^{s,x}-\bar Z^{s,x}_t \;dW_t,\qquad t\in [s,T],\\
  Y_T=\phi(X_T^{s,x}),\\
 Y^{s,x}_t\geq h(X_t^{s,x}),
\end{array}\right.
\end{equation}
which is ``almost'' a reflected BSDE, what is lacking is the requirement that $\bar K^{s,x}$
is the minimal increasing process, namely it is not required the condition
\[
  \int_s^T \left(\bar Y^{s,x}_t-h(X_t^{s,x})\right)\,d\bar K^{s,x}_t=0
\]
Now we have to compare $\bar Y^{s,x}$ with $ Y^{s,x}$. To this aim,
extending a procedure used in \cite{Be}, we compare $ \bar Y^{s,x}$
with the penalized solution $Y^{n,s,x}$ of equation \ref{penaliz}, that we rewrite in integral form,
for $t\in [s,T]$
\begin{equation}\label{penaliz-int}
Y^{n,s,x}_t=\phi (X_T^{s,x})+
\int_t^T\psi(\tau,X^{s,x}_\tau,Y^{n,s,x}_\tau,Z^{n,s,x}_\tau)\;d\tau
+\int_t^T n(Y^{n,s,x}_\tau-h(X^{s,x}_\tau))^-d\tau-\int_t^T Z^{n,s,x}_\tau \;dW_\tau.
\end{equation}
Applying It\^o formula to the process $e^{n(T-t)}Y^{n,s,x}_t$ we get
\begin{equation}\label{penaliz-exp}\left\{\begin{array}{l}
- \dis de^{n(T-t)}Y^{n,s,x}_t=e^{n(T-t)}\psi(t,X^{s,x}_t,Y^{n,s,x}_t,Z^{n,s,x}_t)\;dt
+n e^{n(T-t)}Y^{n,s,x}_t\vee h(X^{s,x}_t)\,dt\\
\qquad\qquad\qquad -e^{n(T-t) }Z^{n,s,x}_t \;dW_t,\qquad t\in [s,T],\\
Y^{n,s,x}_T=\phi(X_T^{s,x})
\end{array}\right.
\end{equation}
Applying It\^o formula to the process $e^{n(T-t)}\bar Y^{s,x}_t$ we get
\begin{equation}\label{bar-exp}\left\{\begin{array}{l}
- \dis de^{n(T-t)}\bar Y^{s,x}_t=n e^{n(T-t)}\bar Y^{s,x}_t+
e^{n(T-t)}\psi(t,X^{s,x}_t,\bar Y^{s,x}_t,\bar Z^{s,x}_t)\;dt
+e^{n(T-t)}\,d\bar K_t^{s,x}\\
\qquad\qquad \qquad-e^{
n(T-t) }\bar Z^{s,x}_t \;dW_t,\qquad t\in [s,T],\\
\bar Y^{s,x}_T=\phi(X_T^{s,x})
\end{array}\right.
\end{equation}
Notice that in (\ref{bar-exp}) we can replace $\bar Y^{s,x}_t$ by $\bar Y^{s,x}_t\vee h(X_t^{s,x})$ (recall that since $\bar u$ is a supersolution to the
obstacle problem (\ref{obstacle}) it holds
$\bar u\geq h$).
Assume for a moment the following lemma.
\begin{lemma}\label{lemma-comparison}
 Let $f^i:\Omega\times[0,T]\times\R\times\Xi\rightarrow \R,\;i=1,2$ satisfy
hypothesis \ref{ip-f-YZ} with $p=2$, fix  $\xi\in L^2_{\calf_T}(\Omega)$  let
$K$ be a progressively measurable nondecreasing processes with $\E\, K_T^2 <\infty$. If
let $(Y^1,Z^1)$ and $(Y^2,Z^2)$ with $Y^i \in  L^2_\calp(\Omega, C([0,T]))$ and $Z^i\in L^2_\calp(\Omega\times[0,T],\Xi)$, $ i=1,2$,
  are the solutions to the following equations of backward
type:
\begin{equation}\label{BSDE-lemma-1}\left\{\begin{array}{l}
  -\dis dY^1_t=f^1(t,Y^1_t,Z^1_t)\;dt+dK_t- Z^1_t \;dW_t,\qquad t\in [0,T],\\
  Y^1_T=\xi,
\end{array}\right.
\end{equation}
\begin{equation}\label{BSDE-lemma-2}\left\{\begin{array}{l}
  -\dis dY^2_t=f^2(t,Y^2_t,Z^2_t)\;dt- Z^2_t \;dW_t,\qquad t\in [0,T],\\
  Y^2_T=\xi.
\end{array}\right.
\end{equation}
and
\begin{equation}\label{delta2}
 \delta_2 f_t: =f^1(t,Y^2_t,Z^2_t)-f^2(t,Y^2_t,Z^2_t)\geq 0,\;d\P\times dt \quad\text{a.s},
\end{equation}
then we have that $Y^1_t\geq Y^2_t$ $\P$-almost surely for any $t\in [0,T]$.
\end{lemma}
By applying lemma \ref{lemma-comparison} to the BSDEs
\ref{penaliz-exp} and \ref{bar-exp} we get a comparison for the processes
$\left(e^{n(T-t)} Y^{n,s,x}_t\right)_{t\in[s,T]}$ and
$\left(e^{n(T-t)}\bar Y^{s,x}_t\right)_{t\in[s,T]}$,
namely we get
\begin{equation}
 e^{n(T-t)}\bar Y^{s,x}_t\geq e^{n(T-t)} Y^{n,s,x}_t
\end{equation}
almost surely and for any time $t$, and consequently
\begin{equation}\label{comparison-Y}
\bar Y^{s,x}_t\geq  Y^{n,s,x}_t.
\end{equation}
Now we let $n\rightarrow\infty$: by \cite{kakapapequ}, section 6,
$Y^{n,s,x}_t\uparrow Y^{s,x}_t$ for any $s\leq t\leq T$ and $\P$-a.s..
So taking $s=t$ in (\ref{comparison-Y}) we finally get
\begin{equation}\label{comparison-u}
\bar u(s,x)\geq  u(s,x),
\end{equation}
for any $\bar u$ supersolution for the obstacle problem \ref{obstacle}.
So the minimality of $u$ is proved: the unique solution to the obstacle problem
\ref{obstacle} is given by formula (\ref{def-di-u}) and the other properties
follows by estimates (\ref{stimaYZ_p}), which passes to the limit as $n\rightarrow\infty$
as stated in proposition \ref{prop-p-RBSDE} on the solution of the RBSDE in terms
of the growth of $\psi$, $h$ and $\phi$.

\qed

In order to complete the proof of theorem \ref{teo-solmild},
we have to prove lemma \ref{lemma-comparison}.

\medskip

\textbf{Proof of Lemma \ref{lemma-comparison}}. We adequate the proof of the classical
comparison theorem for BSDEs given in \cite{kapequ}, Theorem 2.2, to the equations
\ref{BSDE-lemma-1} and \ref{BSDE-lemma-2}. By denoting
\[
\Delta_y f^1_t=\dfrac{f^1(t,Y^1_t,Z^1_t)-f^1(t,Y^2_t,Z^1_t)}{Y^1_t-Y^2_t}\qquad\text{
if }Y^1_t-Y^2_t\neq 0, \qquad
\Delta_y f^1_t=0 \text{ otherwise},
\]
\[
\Delta_z f^1_t=\dfrac{f^1(t,Y^2_t,Z^1_t)-f^1(t,Y^2_t,Z^2_t)}{\vert Z^1_t-Z^2_t\vert^2}
(Z^1_t-Z^2_t)\qquad
\text{ if } Z^1_t-Z^2_t\neq 0,\qquad
\Delta_z f^1_t=0\text{ otherwise},\]

$\delta_2$ as defined in (\ref{delta2}),
$\delta Y_t=Y^1_t-Y^2_t$ and $\delta Z_t=Z^1_t-Z^2_t$
we get
\begin{equation}\label{BSDE-lemma-difference}\left\lbrace\begin{array}{l}
  -\dis d\delta Y_t=\Delta_y f^1_t\,\delta Y_t\;dt
 +(\Delta_z f^1_t)^*\,\delta Z_t\;dt+\delta_2\,f_t\;dt+dK_t- \delta Z_t \;dW_t,\qquad t\in [0,T],\\
  \delta Y_T=0
\end{array}\right.
\end{equation}
We notice that $\Delta_y f^1_t$ and $\Delta_z f^1_t $ are bounded 
and that $\delta_2 f \in  L^2_\calp(\Omega\times[0,T],\R)$.

Multiplying $\delta Y_t$ by $\exp(\int_0^t \Delta_y f^1_{\tau} d\tau)$ and then applying Girsanov theorem we obtain:

\begin{equation}\label{deltaY}
\delta Y_t=\E\left(\rho_{t,T} \left[ \int_t^T \exp\left( \int_t^s \Delta_y f^1_{\tau} d\tau\right) dK_{s }+  \int_t^T\exp\left(\int_t^s \Delta_y f^1_{\tau} d\tau \right)  \delta_2 f_{s} d{s}\right]\right)\end{equation}
where $\rho_{t,T}$ is the Girsanov density:
$$  \rho_{t,T} = \exp \left(\displaystyle \int_t^T (\Delta_z f^1_s)^* dW_s-\frac{1}{2} \int_t^T |\Delta_z f^1_s|^2 ds\right).$$
The claim obviously follows from (\ref{deltaY}) being $(K)$ non decreasing and $\delta_2 Y$ non negative. \qed

\section{The Optimal Control-Stopping problem} \label{relfond}
 An \textit{Admissible Control System} is a set
 $$\cals=(  \Omega^{\cals},\mathcal{F}^{\cals}, (\mathcal{F}^{\cals}_t)_{t\geq 0},\mathbb{P}^{\cals},( W_t^{\cals})_{t
\geq 0})$$
where $(  \Omega^{\cals},\mathcal{F}^{\cals}, (\mathcal{F}^{\cals}_t)_{t\geq 0},\mathbb{P}^{\cals})$
 is a complete probability space endowed with a filtration satisfying the usual assumptions and $( W_t^{\cals})_{t
\geq 0})$ is a cylindrical Wiener process in  $\Xi$.
Fixed a closed subset $U$ of
a normed space $U_0$ an \textit{admissible control}  in the setting $\cals$ is any $\left(  \mathcal{F}^{\cals}_{t}\right)  $-predictable
process  $\alpha:\Omega^{\cals}\times [0,T]\rightarrow U$. The set of all admissible controls will be denoted by $\calu^{\cals}$.

We fix a function $R: H\times U \rightarrow \Xi$ bounded, continuous such that:
\begin{equation}\label{lip_di_R}
|R(\alpha,x)-R(\alpha,x')| \leq |x-x'| \qquad \forall u\in U,\; x,x' \in H
\end{equation}

Given an admissible setting $\cals$ and an admissible control  $\alpha\in \calu^{\cals}$ and fixed $x\in H$, $s\in [0,T]$ by $X^{\alpha,s,x}$ we will denote the solution
to the following stochastic differential equation in a Hilbert space $H$%
\begin{equation}
\left\{
\begin{array}
[c]{l}%
dX_{t}^{\alpha, s,x}=AX_t^{\alpha,s,x} dt +F(t,X_t^{\alpha,s,x})dt+G(t,X_t^{\alpha,s,x})(R(X^{\alpha, s,x}_t,\alpha_t) dt +dW^{\cals}_t),\text{ \ \ \ }t\in\left[  s,T\right] \\
X_{s}^{\alpha, s,x}=x\in H.
\end{array}
\right.  \label{SDEcontrol}%
\end{equation}
Moreover given $l: [0,T]\times H \times U_0\rightarrow \R$ we introduce the cost
 functional:
\begin{equation}
J\left(  s,x,\tau,\alpha \right)  =\E\int_{s}%
^\tau l\left(  r,X_{r}^{\alpha,s,x},\alpha_{r}\right)  dr+
\E[\phi\left(  X_{T}^{\alpha,s,x} \right)\chi_{\left\lbrace \tau =T\right\rbrace} ] +
\E [h\left(\tau , X_{\tau}^{\alpha,s,x}\right) \chi_{\left\lbrace\tau <T\right\rbrace}] ,
 \label{costo}%
\end{equation}
that we wish to maximize over all admissible control  $\alpha\in \calu^{\cals}$
and over all $ \{\calfì^{\cals}_t\}_t$-stopping times $\tau$ satisfying
$t\leq \tau \leq T$.

For $s\in [0,T]$, $x\in H$, $z\in \Xi^*$ we define the
hamiltonian function in the usual way as
\begin{equation}\label{defhamiltonian}
\psi(s,x,z)=\sup_{\alpha\in\calu} \{zR(x,\alpha) +l(s,x ,\alpha) \}.
\end{equation}
We notice that since $R$ is bounded $\psi$ is Lipschitz with respect to $z$. 
$ $
We will assume throughout this section that $A$, $F$ and $G$ verify Hypothesis (\ref{ip-forward}) and that $\phi$, $\psi$ and $h$ verify  Hypothesis (\ref{ip su psi}). 
Moreover we assume that $|l(s,x,\alpha)|\leq c(1+|x|^r)$ for some $c,r>0$.

We notice that under the above assumptions, fixed $s\in [0,T]$ and $x\in H$ then for all $\alpha\in \calu^{\cals}$ there exists a unique mild solution  $X^{\alpha,s,x}$ to equation (\ref{SDEcontrol}). Moreover  $X^{\alpha,s,x}\in L^p_\calp(\Omega,C([s,T],H))$ for all $p \geq 1$, see \cite{FuTe}. Consequently
$J(s,x,\tau,\alpha)$ is a well defined real number for all $\alpha\in \calu^{\cals}$
and all $ \{\calfì^{\cals}_t\}_t$-stopping time $\tau \leq T$.
We also notice that $X^{\alpha,s,x}$ is adapted to the filtration generated by $(W^{S}_t)$.

By the Girsanov theorem, there exists a
probability measure $\P^{\cals,\alpha}$  such that  the process $$W^{\cals,\alpha}_{t}:=W^{\cals}_t+\int_s^tR(X^{\alpha,s,x}_r, \alpha_r) dr \quad t\geq s $$
is a cylindrical $\P^{\cals,\alpha}$-Wiener process in $\Xi$. We denote by 
$(  \mathcal{F}^{\cals,\alpha}_t)  _{t\geq
s}$ its natural filtration, augmented in the usual way.  $X^{\alpha,s,x}$ satisfies the following equation:
\begin{equation}
\left\{
\begin{array}
[c]{l}%
dX_{t}^{\alpha, s,x}=AX_t^{\alpha,s,x} dt +F(t,X_t^{\alpha,s,x})dt+G(t,X_t^{\alpha,s,x})dW^{\cals,\alpha} _t,\text{ \ \ \ }t\in\left[  s,T\right] \\
 X_{s}^{\alpha, s,x}=x.
\end{array}
\right.  \label{SDE_Girsanov}%
\end{equation}
Consequently (notice that the above equation enjoys
strong existence, in probabilistic sense, and pathwise  uniqueness) $X_{t}^{\alpha, s,x}$ turns out to be adapted to $(  \mathcal{F}^{\cals,\alpha}_t)_{t\geq
s}.$

In $\left(  \Omega^{\cals},\mathcal{F}^{\cals},
(  \mathcal{F}_t^{\cals,\alpha})_{t\geq
0},\mathbb{P}^{\cals,\alpha}\right)$ we consider the solution $(\widetilde Y^{s,x}, \widetilde Z^{s,x},\widetilde K^{s,x})$ of the following reflected backward stochastic
differential equation:
\begin{equation}\label{RBSDE_Girsanov}\left\{\begin{array}{l}
  -\dis d\widetilde Y_t^{s,x}=\psi(s,X_t^{\alpha,s,x},\widetilde Z_t^{s,x})\;dt
+d\widetilde K_t^{s,x}-\widetilde Z_t^{s,x} \;dW^{\cals,\alpha}_t,\qquad t\in [0,T],\\
 \widetilde Y_T^{s,x}=\phi(X_T^{\alpha,s,x}),\\
 \widetilde Y_t^{s,x}\geq h(t,X_t^{\alpha,s,x}),\\
 \int_0^T (\widetilde Y_t^{s,x}-h(t,X_t^{\alpha,s,x}))d\widetilde K_t^{s,x}=0,
\end{array}\right.
\end{equation}
We omit to indicate the dependence on the admissible setting $\cals$ and on the admissible control $\alpha$ since the 
 law of $\widetilde Y^{s,x}$, $\widetilde Z^{s,x}$ and $\widetilde K^{s,x}$
is uniquely determined by $A$, $F$, $G$, $x$, $\psi$ and $\phi$,
and does not depend on the probability space and on the Wiener process,
and in particular $\widetilde Y^{s,x}_s$ is a real number that does not depend on $\cals$ and on $\alpha$.
We argue as in \cite{kakapapequ}, proposition 2.3. Rewriting 
 (\ref{RBSDE_Girsanov}) in terms of the original noise $(W^{\cals})$ and integrating it between $s$ and any  $(\calf^{\cals}_s)$-stopping time $\tau$, we get that $\P^{\cals,\alpha}$-a.s., and consequently
$\P^{\cals,\alpha}$-a.s.,
$$
\widetilde Y_s^{s,x} = \widetilde Y_\tau^{s,x}+
\int_s^\tau\psi(r,X_r^{\alpha,s,x},\widetilde Z_r^{s,x})\;dr
+\widetilde K_\tau^{s,x}-\widetilde K_t^{s,x}
-\int_s^\tau\widetilde Z_r^{s,x}\;dW^{\cals}_r-\int_s^\tau\widetilde Z_r^{s,x}\;R(X_r^{\alpha,s,x},\alpha_r) dr \nonumber$$
Noticing that
 $(\int_0^t \widetilde Z_r^{\alpha,s,x} \;dW^{\cals}_r)_{t\geq 0}$ is a
$\P^{\cals}$-martingale and that  
$\widetilde Y_r^{\alpha,s,x}\geq h(r, X^{\alpha,s,x}_r)$ by computing expectation with respect to  $\P^{\cals}$ we get:
\begin{align*}
\widetilde Y_s^{s,x} & \geq
 \E  \int_s^\tau\psi(r,X_r^{\alpha,s,x},\widetilde Z_r^{\alpha,s,x})\;dr
-\E\int_s^\tau \widetilde Z_r^{\alpha,s,x}\alpha_r\;dr \\ \nonumber
&+\E[\widetilde K_\tau^{s,x}-\widetilde K_t^{s,x}]+\E h(\tau,X_\tau^{\alpha,s,x})\chi_{\{\tau<T\}}+ \E \phi(X_T^{\alpha,s,x})\chi_{\{\tau=T\}}],\\ \nonumber
\end{align*}
Finally adding and subtracting the current cost we have:
\begin{equation}\label{relfond1}
\widetilde Y_s^{s,x}  \geq
 J(s,x,\tau,\alpha)+\E  \int_s^\tau\left[\psi(r,X_r^{\alpha,s,x},\widetilde Z_r^{\alpha,s,x})- l(r,X_r^{\alpha,s,x},u_r)- Z_r^{s,x}u_r\right]\;dr 
 +\E[\widetilde K_\tau^{s,x}-\widetilde K_s^{s,x}]. \\
 \end{equation}

We have therefore proved the following result
\begin{theorem}\label{teo_rel_fond} For every admissible setting $\cals$ and every admissible control $u\in \calu^{\cals}$ we have:
$$ 
 J(s,x,\tau,\alpha)\leq  \widetilde Y_s^{s,x} $$
 moreover the equality holds if and only if
 \begin{eqnarray}\label{cond-suff}
& \psi(r,X_r^{\alpha,s,x},\widetilde Z_r^{s,x})- l(r,X_r^{\alpha,s,x},u_r)- Z_r^{s,x}\alpha_r=0,\quad \P-\hbox{a.s. for a.e.} r\in [s,\tau]\\
&  K_\tau^{s,x}-\widetilde K_s^{s,x}=0
  ,\quad \P-\hbox{a.s.} \\
  & \widetilde Y_{\tau}^{s,x}I_{\{\tau<T\}}=h(\tau,X_{\tau}^{\alpha,s,x})I_{\{\tau<T\},}\quad \P-\hbox{a.s.}\label{cond-suff-3}
  \end{eqnarray}
  \end{theorem}
  \begin{remark} \label{corollario_tempo_ottimo} Fixed an admissible setting $\cals$ and ad admissible control
  $\alpha\in \calu^{\cals}$ let $\bar{\tau}$ be define as 
  Let us consider
\begin{equation}\label{deftempo_ottimo}
\overline{\tau}=\inf\{ t\leq r\leq T: \widetilde Y_r^{s,x}=h(r,X_r^{\alpha,s,x})\}\wedge T.
\end{equation}
The condition
$ \int_0^T (\widetilde Y_t^{s,x}-h(t,X_t^{\alpha,s,x}))d\widetilde K_t^{s,x}=0$
together with continuity and monotonicity of $\widetilde K$ imply that
$$
\widetilde K_{\overline{\tau}}^{s,x}-\widetilde K_t^{s,x}=0.
$$
Moreover (\ref{cond-suff-3}) follows by definition. Consequently we have:
\begin{equation}\label{relfond2}
\widetilde Y_s^{s,x}  =
 J(s,x,\bar{\tau},\alpha)+\E  \int_s^{\bar{\tau}}\left[\psi(r,X_r^{\alpha,s,x},\widetilde Z_r^{\alpha,s,x})- l(r,X_r^{\alpha,s,x},\alpha_r)- Z_r^{\alpha,s,x}\alpha_r\right]\;dr.
\end{equation}
   \end{remark}
   Taking into account equations (\ref{SDE_Girsanov}), (\ref{RBSDE_Girsanov}) and Proposition \ref{teo-solmild} the above results can be refomulated as follows
   \begin{corollary}\label{cor_rel_fond}
 Let $u$ be the minimal mild supersolution to the obstacle problem and let $\zeta$ be any element of its generalized gradient. Given any admissible setting $\cals$ and any admissible control $\alpha\in \calu^{\cals}$ we have:
$$ 
 J(s,x,\tau,\alpha)\leq u(s,x) $$
 moreover the equality holds if and only if
 \begin{eqnarray}\label{cond-suff-HJB}
& \psi(r,X_r^{\alpha,s,x}, \zeta (r,X_r^{\alpha,s,x}))- l(r,X_r^{u,s,x},\alpha_r)- \zeta (r,X_r^{\alpha,s,x})\alpha_r=0,\quad \P-\hbox{a.s. for a.e.} r\in [s,\tau] \nonumber\\
&  K_\tau^{s,x}-\widetilde K_s^{s,x}=0
  ,\quad \P-\hbox{a.s.},\nonumber\\
  & u(\tau, X_{\tau}^{\alpha,s,x})I_{\{\tau<T\}}=h(\tau,X_{\tau}^{\alpha,s,x})I_{\{\tau<T\},}\quad \P-\hbox{a.s.}.\nonumber
  \end{eqnarray}
  Finally if 
  \begin{equation}\label{deftempo_ottimo-HJB}
\overline{\tau}=\inf\{ t\leq r\leq T:  u(r,X_r^{\alpha,s,x})=h(r,X_r^{\alpha,s,x})\}\wedge T.
\end{equation}
then the equality holds if and only if
(\ref{cond-suff-HJB}) holds.
  
\end{corollary}
  We come now to the existence of optimal controls. We shall exploit the weak formulation of the control problem and select a  suitable admissible setting $\bar{\cals}$. We assume the following
  \begin{hypothesis} \label{hyp_selezione}
  The minimum in the definition (\ref{defhamiltonian}) is attained for all $t\in [s,T]$, $x\in H$ and $z\in \Xi^*$  e.g. if we define
\begin{equation}\label{defdigammagrande}
\Gamma(s,x,z)=\left\{ \alpha\in \calu: zR(x,\alpha)+l(s,x ,\alpha)= \psi(s,x,z)\right\}
\end{equation}
 then
$\Gamma(s,x,z) \neq \emptyset$ for every $s\in [0,T]$, every $x\in
H$ and every $z\in \Xi^*$.
\end{hypothesis}
\begin{remark}
 By \cite{AubFr}, see Theorems 8.2.10 and
8.2.11, under the above assumption $\Gamma$ always admits a measurable selection, i.e. there exists
a measurable function $\gamma: [0,T]\times H \times \Xi^*
\rightarrow U$ with $\gamma(s,x,z)\in \Gamma(s,x,z)$ for every
$s\in [0,T]$, every $x\in H$ and every $z\in \Xi^*$.

Moreover we notice that if $U$ is compact then Hypothesis \ref{hyp_selezione} always hold
\end{remark}
\begin{theorem}
Assume Hypothesis \ref{hyp_selezione} and fix a measurable selection $\gamma$ of $\Gamma$, $s\in [0,T]$, $x\in H$ and an element $\zeta $ of the generalized gradient of the minimal supersolution $u$ of the obstacle problem (\ref{obstacle}); then there exists at least an admissible setting $\bar{\cals}$ in which the closed loop equation
\begin{equation}
\left\{
\begin{array}
[c]{l}%
d\bar{X}_{t}=A\bar{X}_tdt +F(t,\bar{X}_t)dt+G(t,\bar{X}_t)[R(t,\gamma(t,\zeta (t, \bar{X}_t)))+dW^{\cals}_t],\text{ \ \ \ }t\in\left[  s,T\right] \\
 \bar X_{s}=x.
\end{array}
\right.  \label{closed_loop}%
\end{equation}
admits a mild solution.
\end{theorem}

\textbf{Proof.} We fix any admissible setting  $$\cals=(  \Omega^{\cals},\mathcal{F}^{\cals}, (\mathcal{F}^{\cals}_t)_{t\geq 0},\mathbb{P}^{\cals},( W_t^{\bar\cals})_{t
\geq 0})$$ and consider the uncontrolled forward SDE
   \begin{equation}
\left\{
\begin{array}
[c]{l}%
dX_{t}=AX_t dt +F(t,X_t)dt+G(t,X_t)dW^{\cals}_t,\text{ \ \ \ }t\in\left[  s,T\right] \\
 X_{s}=x.
\end{array}
\right.  \label{SDE_senza_controllo}%
\end{equation}

By the Girsanov theorem, there exists a
probability measure $\hat{\P} $  such that  the process $$\hat{W}_{t}:=W^{\cals}_t-\int_s^tR(X^{s,x}_r, \zeta(s,X^{s,x}_r)) dr 
\quad t\geq s $$
is a cylindrical $\hat{\P}$-Wiener process in $\Xi$. We denote by 
$( \hat{ \mathcal{F}}_t)  _{t\geq
s}$ its natural filtration, augmented in the usual way. Clearly 
 ${X}$ solves 
\begin{equation}
\left\{
\begin{array}
[c]{l}%
\displaystyle d{X}_{t}=A{X}_{t} dt +F(t,{X}_{t})dt+G(t,{X}_{t})[R({X}_t, \gamma(t,\zeta (t, {X}_t)dt+d\hat{W} _t],\text{ \ \ \ }t\in\left[  s,T\right] \\
\displaystyle\hat{ X}_{s}=x.
\end{array}
\right.%
\end{equation}
and  $(  \Omega^{\cals},\mathcal{F}^{\cals}, (\hat{\mathcal{F}}_t)_{t\geq 0},\hat{\mathbb{P}},( \hat{W}_t)_{t
\geq 0})$ is the desired admissible system.
\qed

We finally get the following
\begin{theorem}
Assume Hypothesis \ref{hyp_selezione} and fix a measurable selection $\gamma$ of $\Gamma$, $s\in [0,T]$, $x\in H$ and an element $\zeta $ of the generalized gradient of the minimal supersolution $u$ of the obstacle problem (\ref{obstacle}). Moreover let  $\bar{\cals}$ be an admissible setting in which the closed loop equation (\ref{closed_loop}) admits a mild solution then there exists $\bar{\alpha}\in \calu^{\bar{\cals}}$
 and an $(\calf^{\bar{\cals}})$ stopping time $\bar{\tau}$ for which
$$
\widetilde Y_s^{s,x}=u(s,x)=
 J(s,x,\bar{\tau},\bar{\alpha}).
$$
\end{theorem}

\textbf{Proof:} Just let $\bar{X}$ be the mild solution of equation (\ref{closed_loop})and define $\bar{\alpha}=\gamma(t, \zeta (t,\bar{X}_t))$ clearly $\bar{X}_t= X^{\bar{\alpha},s,x}$ and relation (\ref{cond-suff-HJB}) holds. Thus by Corollary \ref{cor_rel_fond} it is enough to choose 
$$
\overline{\tau}=\inf\{ t\leq r\leq T:  u(r,\bar{X}_r)=h(r,\bar{X}_r)\}\wedge T. $$
$ $ \qed



\begin{thebibliography}{11}

 \bibitem{AubFr} {J.P. Aubin, H. Frankowska},
    \textit{ Set-valued analysis},
      {Systems \& Control: Foundations \& Applications},
Vol. {2}, {Birkh\"auser Boston Inc.}, {Boston, MA},     {1990},

\bibitem{Be} A. Bensoussan, \textit{Stochastic control by functional analysis methods.}
 Studies in Mathematics and its Applications, 11.
North-Holland Publishing Co., Amsterdam-New York, 1982.


%


 \bibitem {DP1}G. Da Prato, J. Zabczyk,\textit{ Stochastic equations in
 infinite dimensions. }Encyclopedia of Mathematics and its Applications 44,
 Cambridge University Press, 1992.

 \bibitem {DP3} G. Da Prato, J. Zabczyk,
\textit{ Second order partial differential equations in Hilbert spaces. L}. ondon Mathematical Society Lecture Note Series, 293. Cambridge University Press, Cambridge, 2002

\bibitem {DelMey} C. Dellacherie, P. A. Meyer \textit{Probability and Potential B: Theory of Martingales }, North-Holland Amsterdam (1982).

\bibitem{kakapapequ} N. El Karoui, C. Kapoudjian, E. Pardoux,  S. Peng, M. C.
Reflected solutions of backward SDE's, and related obstacle problems for PDE's.
Ann. Probab.  \textbf{25}  (1997),  no. 2, 702--737.

\bibitem{kapequ} N. El Karoui, S. Peng, M. C. Quenez, Backward
stochastic differential equations in finance. Mathematical Finance
\textbf{7} (1997), 1-71.




\bibitem{FuTe} M. Fuhrman, G. Tessitore, Nonlinear Kolmogorov
equations in infinite dimensional spaces: the backward stochastic
differential equations approach and applications to optimal
control. Ann. Probab. \textbf{30} (2002), 1397--1465.

\bibitem{FuTe_Bismut} M. Fuhrman, G. Tessitore, The Bismut-Elworthy formula for backward SDEs and applications to nonlinear Kolmogorov equations and control in infinite dimensional spaces. Stoch. Stoch. Rep. 
\textbf{74} (2002), no. 1-2, 429--464.


\bibitem{FuTe-ell} M. Fuhrman, G. Tessitore, Infinite horizon
 backward stochastic differential equations
and elliptic equations in Hilbert spaces. Ann. Probab. \textbf{32}
(2004), 607--660.


\bibitem{FuTeGen} M. Fuhrman, G. Tessitore, Generalized directional gradients, backward stochastic differential equations and mild solutions of semilinear parabolic equations. Appl. Math. Optim. 51 (2005), no. 3, 279--332.

\bibitem{HuTess} Y. Hu, G. Tessitore, BSDE on an infinite horizon
and elliptic PDEs in infinite dimension.  NoDEA Nonlinear
Differential Equations Appl.  14  (2007),  no. 5-6, 825--846


\bibitem{KaSh} I. Karatzas, S.E. Shreve, Steven, \textit{Brownian motion and
stochastic calculus. Second edition.} Graduate Texts in Mathematics, 113.
Springer-Verlag, New York.

\bibitem{kesw} D. Kelome, A. Swiech, Viscosity solutions of an infinite-dimensional             Black-Scholes-Barenblatt equation,{Appl. Math. Optim.}, \textbf{47} (2003),253--278.


\bibitem{Mas} F. Masiero, Semilinear Kolmogorov equations and
applications to stochastic optimal control, Appl. Math. Optim.,
51 (2005), pp.~201--250.

\bibitem{Mas1} F. Masiero, Infinite horizon stochastic optimal control problems with degenerate noise and elliptic equations in Hilbert spaces. Appl. Math. Optim. 55 (2007), no. 3, 285-326





\bibitem{PaPe-90} E. Pardoux, S. Peng,
Adapted solution of a backward stochastic differential equation.
Systems and Control Lett. {\bf 14}, 1990, 55-61.

\bibitem{PaPe} E. Pardoux, S. Peng,
Backward stochastic differential equations and quasilinear
parabolic partial differential equations, in: \textit{ Stochastic
partial differential equations and their applications}, eds. B.L.
Rozowskii, R.B. Sowers, 200-217, Lecture Notes in Control Inf.
Sci. 176, Springer, 1992.

\bibitem{YongZhou} {J. Yong, X. Y. Zhou},
     \textit{Stochastic controls, {Hamiltonian systems and HJB equations}},
    Applications of Mathematics,
   Springer, {New York},
      (1999).


\end{thebibliography}
\end{document}